\documentclass[12pt,twoside]{article}
\usepackage{times}
\usepackage{amsmath,amssymb}
\usepackage{color}

\allowdisplaybreaks

\def\p{\partial}
\def\ve{\varepsilon}
\def\f{\frac}
\def\na{\nabla}
\def\la{\lambda}

\def\t{\tilde}

\def\O{\Omega}
\def\th{\theta}
\def\g{\gamma}
\def\G{\Gamma}
\def\si{\sigma}
\def\dl{\delta}
\def\p{\partial}

\def\ve{\varepsilon}
\def\f{\frac}

\def\Dl{\Delta}

\def\na{\nabla}
\def\la{\lambda}

\def\t{\tilde}

\def\th{\theta}

\def\g{\gamma}
\def\G{\Gamma}

\def\si{\sigma}
\def\dl{\delta}

\def\ds{\displaystyle}

\def\t{\tilde}

\def\th{\theta}

\def\g{\gamma}
\def\G{\Gamma}
\def\si{\sigma}

\def\dl{\delta}
\def\Dl{\Delta}
\def\i{\infty}

\def\Div{\text{div }}
\def\re{\text{Re }}

 \allowdisplaybreaks

\begin{document}
 \footskip=0pt
 \footnotesep=2pt
\let\oldsection\section
\renewcommand\section{\setcounter{equation}{0}\oldsection}
\renewcommand\thesection{\arabic{section}}
\renewcommand\theequation{\thesection.\arabic{equation}}
\newtheorem{claim}{\noindent Claim}[section]
\newtheorem{theorem}{\noindent Theorem}[section]
\newtheorem{lemma}{\noindent Lemma}[section]
\newtheorem{proposition}{\noindent Proposition}[section]
\newtheorem{definition}{\noindent Definition}[section]
\newtheorem{remark}{\noindent Remark}[section]
\newtheorem{corollary}{\noindent Corollary}[section]
\newtheorem{example}{\noindent Example}[section]

\title{Large time asymptotic behavior of the compressible Navier-Stokes Equations in partial
Space-Periodic Domains}

\author{Yin Huicheng$^{*}$, \quad Zhang Lin$^{*}$,\quad Zhu Lu
\footnote{This project was supported by the NSFC (
No.~11025105), and by the
Priority Academic Program Development of Jiangsu Higher Education
Institutions.}\vspace{0.5cm}\\
\small  (Department of Mathematics and
IMS, Nanjing University, Nanjing 210093, China.)\\
\vspace{0.5cm}
}

\date{}
\maketitle

\centerline {\bf Abstract} \vskip 0.3 true cm

In this paper, we study the large time behavior of the 3-D isentropic compressible Navier-Stokes equation in the partial
space-periodic domains, and simultaneously show that the related profile systems can be described by like
Navier-Stokes equations with suitable ``pressure'' functions in lower dimensions. Our proofs are based on the energy
methods together with some delicate analysis on the corresponding linearized problems.

\vskip 0.3 true cm

\vskip 0.3 true cm


\vskip 0.4 true cm
\centerline{\bf $\S 1$. Introduction  and main results}
\vskip 0.3 true cm

In this paper, we consider the 3-D  isentropic compressible Navier-Stokes equation for
$(t, z)\in [0, +\infty)\times\O$:
\begin{equation}
\left\{
\begin{aligned}
&\partial_{t}\rho+\Div m=0,\\
&\partial_{t}m+\Div(\ds\f{m\otimes m}{\rho})+\nabla p(\rho)=\mu\Delta(\ds\f{m}{\rho})+(\mu+\mu')\nabla\Div(\ds\f{m}{\rho}),
\end{aligned}
\right.\tag{1.1}
\end{equation}
where $\Omega=\Bbb T^\ell\times\Bbb R^{3-\ell}$, $\Bbb T^\ell=[0,2\pi]^\ell$ is the $\ell$-dimensional torus $(1\le\ell\le 3)$, $z=(z_1, z_2, z_3)$, $\rho=\rho(t,z)$ is the density, $m=(m^{1},m^{2},m^{3})(t,z)$ is the  momentum, $p=p(\rho)$ is the pressure, $\mu$ and $\mu'$ are the first and the second viscosity coefficient respectively, which satisfy $\mu>0$ and $\ds\f{2}{3}\mu+\mu'\ge0$. The initial data of (1.1)
are given as follows
$$
(\rho, m)(0,z)=(1+\rho_0(z), m_0(z)),\eqno{(1.2)}
$$
where $(\rho_0, m_0)\in H^4(\Omega)\times(H^4(\Omega))^3$, and $1+\rho_0(z)>0$ for $z\in\O$.

It is obvious that $(\rho, m)=(1,0)$ is a steady solution of (1.1) with the initial data $(\rho, m)(0,z)$ $=(1, 0)$.
We will be concerned with the perturbation problem of (1.1) to this constant state. Denote by
$$
\nu_1=\mu,\ \nu_2=\mu+\mu',\ \g=\sqrt{p'(1)}.
$$
As in [14, 15], if we set $\phi=\g(\rho-1)$, then (1.1)-(1.2) can be rewritten as
\begin{equation}
\left\{
\begin{aligned}
&\partial_{t}\phi+\g\Div m=0,\\
&\p_t m-\nu_1\Dl m-\nu_2\na\Div m+\g\na\phi=G(\phi,m),\\
&(\phi,m)(0,z)=(\phi_0,m_0)(z),
\end{aligned}
\right.\tag{1.3}
\end{equation}
where $\phi_0(z)=\g\rho_0(z)$ and
\begin{align*}
G(\phi,m)&=-\Div\biggl(\f{\g}{\phi+\g}m\otimes m\bigg)-\nu_1\Dl\bigg(\f{\phi}{\phi+\g}m\bigg)-\nu_2\na\Div\bigg(\f{\phi}{\phi+\g}m\bigg)\\
&\quad -\na\bigg\{\f{\phi^2}{\g^2}\int_0^1(1-\th)^2p''\big(1+\f{\th\phi}{\g}\big)d\th\bigg\}.
\end{align*}

We now introduce some notations for later uses. The Fourier transformation  of the function $f\in L^1(\Bbb T^\ell\times\Bbb R^{3-\ell})$ is denoted by
$$
\mathcal {F}(f)(k,\xi)=\hat f(k,\xi)=\int_{\Bbb T^\ell}\int_{\Bbb R^{3-\ell}}e^{-i(k\cdot x+\xi\cdot y)}f(x,y)dxdy,
$$
where $z=(x, y)$, $x=(z_1,\cdots,z_l) \in \Bbb T^\ell,\ y=(y_1, \cdots, y_{3-\ell})=(z_{\ell+1},\cdots,z_3)\in\Bbb R^{3-\ell},\ k=(k_1,\cdots,k_\ell)\in\Bbb Z^\ell,\ \xi=(\xi_1,\cdots,\xi_{3-\ell})\in \Bbb R^{3-\ell}$. The inverse Fourier transformation of the sequence $\{g(k,\xi)\}_{k\in \Bbb Z^l}$
is defined as
$$
(\mathcal {F}^{-1}g)(z)=\f{1}{(2\pi)^\ell}\sum_{k\in\Bbb Z^\ell}e^{i k\cdot x}\int_{\Bbb R^{3-\ell}}e^{i\xi\cdot y}g(k,\xi)d\xi.
$$
Write the mean value of $f(z)$ over $\Bbb T^\ell$ as $\bar f(y)$:
$$
\bar f (y) =\f{1}{(2\pi)^\ell}\int_{\Bbb T^\ell}f(x, y)dx.$$
In addition, we define
$$
\Div' v =\p_{y_1} v_1+\cdots+\p_{y_{3-\ell}}v_{3-\ell},\ \nabla'=(\p_{y_1},\cdots,\p_{y_{3-\ell}})^T,\ \Delta'=\p_{y_1}^2+\cdots+\p_{y_{3-\ell}}^2,
$$
where $v=(v_1,\cdots, v_{3-\ell})^T$. And we denote $\|u(t,\cdot)\|_{L^p(\Omega)}$ as $\|u\|_{L^p}$ for $1\le p\le\infty$.

For different $\ell$, our main results in this paper are
\vskip 0.2 true cm
{\bf Theorem 1.1. }{\it For $\ell=1$, if $u_0\in H^4(\Bbb T \times \Bbb R^2)$ and $\|u_0\|_{H^4\cap L^1}\leq\ve$, then for small
$\ve>0$, (1.3) has a global solution $u(t,z)=(\phi,m)(t,z)\in C([0, +\infty), H^4(\Bbb T \times \Bbb R^2))\cap
C^1([0, +\infty), H^2(\Bbb T \times \Bbb R^2))$
satisfying for $t\rightarrow +\infty$
\begin{align*}
&\|\p_{z}^k u\|_{L^2} =O(t^{-\f{1}{2}-\f{k}{2}}),\ \ k=0,1, \tag{1.4}\\
&\|\p_{z}^k u\|_{L^2} =O(t^{-\f{4-k}{3}}),\ \ k=2,3,4, \tag{1.5}\\
&\|u(t,z)-\eta(t,y)\|_{L^2}=O(t^{-1}),\tag{1.6}
\end{align*}
where the profile $\eta(t,y)=(\sigma, w)(t,y)$ satisfies the following system
\begin{equation}
\left\{
\begin{aligned}
&\partial_{t}\sigma +\g \Div 'w'=0,\\
&\p_tw_1-\nu_1\Delta 'w_1=-\Div' ( w_1 w')\\
&\p_tw'-\nu_1\Delta 'w'-\nu_2 \nabla'\Div'w'+\g\nabla'\sigma= -\Div'(w'\otimes w')- \alpha\na'(\sigma^2),\\
&(\sigma,w)(0,y)=(\bar\phi_0,\bar m_0)(y),
\end{aligned}
\right.\tag{1.7}
\end{equation}
here and below $\alpha=\ds\f{p''(1)}{2\gamma^2}>0$, $w'=(w_2, w_3)$.
}

\vskip 0.2 true cm

{\bf Theorem 1.2. }{\it For $\ell=2$,  if $u_0\in H^4(\Bbb T^2 \times \Bbb R)$ and $\|u_0\|_{H^4\cap L^1}\leq\ve$, then for small
$\ve>0$, (1.3) has a global solution $u(t,z)=(\phi, m)(t,z)\in C([0, +\infty), H^4(\Bbb T^2 \times \Bbb R))\cap
C^1([0, +\infty), H^2(\Bbb T^2 \times \Bbb R))$ satisfing for $t\rightarrow +\infty$
\begin{align*}
&\|\p_{z}^k u\|_{L^2}=O(t^{-\f{1}{4}-\f{k}{2}}),\ \ k=0,1, \tag{1.8}\\
&\|\p_{z}^k u\|_{L^2} =O(t^{-\f{4-k}{4}}),\ \ k=2,3,4, \tag{1.9}\\
&\|u(t,z)-\eta(t,y)\|_{L^2}=O(t^{-\f{3}{4}+\delta}), \ \forall\delta >0, \tag{1.10}
\end{align*}
where $\eta(t,y)=(\sigma, w)(t,y)$ satisfies the following equations
\begin{equation}
\left\{
\begin{aligned}
&\partial_{t}\sigma+\g\p_{y}w_3=0,\\
&\p_tw_j-\nu_1\p^2_{y}w_j=-\p_{y}(w_iw_3),\ j=1,2,\\
&\p_tw_3-(\nu_1+\nu_2)\p^2_{y}w_3+\g\p_{y}\sigma=-\p_{y} w_3^2-\alpha \p_{y} \sigma^2,\\
&(\sigma,w)(0,y)=(\bar \phi_0,\bar m_0)(y).
\end{aligned}
\right.\tag{1.11}
\end{equation}
}

\vskip 0.2 true cm

{\bf Theorem 1.3. }{\it For $\ell=3$,  if $u_0\in H^4(\Bbb T^3)$ and $\|u_0\|_{H^4}\leq\ve$, then for small
$\ve>0$, (1.3)
has a global solution $u(t,z)=(\phi,m)(t,z)\in C([0, +\infty), H^4(\Bbb T^3))\cap
C^1([0, +\infty), H^2(\Bbb T^3))$ satisfying $\|u-\bar u_0\|_{H^1}\le Ce^{-a_0 t}\|u_0\|_{H^4}$, where $a_0>0$ is some determined number.}
\vskip 0.2 true cm

{\bf Remark 1.1. }{\it  The local existence and uniqueness of the solution to the problem (1.3) have been shown in [16]. In addition,
the global well-posedness of (1.3) was obtained by the Matsumura-Nishida energy method in [13-15]. So the purpose of this paper
is to study the large-time behavior of the solution to (1.3) and find the corresponding profile systems.}

\vskip 0.2 true cm

{\bf Remark 1.2. }{\it For
the general $N-$dimensional ($N\geq4$) isentropic  compressible Navier-Stokes equation (1.3), by the same
analysis as in Theorem 1.1-Theorem 1.3, we can get the following conclusions: Assuming $u_0=(\phi_0, m_0)\in H^{s}(\Bbb T^\ell\times\Bbb R^{N-\ell})$
and $\|u_0\|_{H^s\cap L^1}\leq\ve$, where $s\geq[\f{N}{2}]+3$, $1\le \ell\le N-1$, and $\ve>0$ is small enough,
then
\begin{align*}
&\|\p_{z}^k u\|_{L^2}=O(t^{-\f{N-l}{4}-\f{k}{2}}),\ \ k=0,1,\\
&\|u(t,z)-\eta(t,y)\|_{L^2}=O(t^{-\f{N-\ell+2}{4}}),
\end{align*}
where $u=(\phi, m)$, and the profile $\eta(t,y)=(\si,w)(t,y)$ satisfies the following partial differential equation system
\begin{equation*}
\left\{
\begin{aligned}
&\partial_{t}\sigma +\g \Div 'V=0,\\
&\p_tU-\nu_1\Delta'U=-\Div'(V\otimes U)\\
&\p_tV-\nu_1\Delta'V-\nu_2\nabla'\Div'V+\g\nabla'\sigma= -\Div'(V\otimes V)-\alpha\na'(\sigma^2),\\
&(\sigma,w)(0,y)=(\bar\phi_0,\bar m_0)(y),\\
\end{aligned}
\right.
\end{equation*}
here $w=(U, V)$, $U=(w_1,\ldots,w_\ell)$ and $V=(w_{\ell+1},\ldots,w_{N})$.

If $\ell=N$, the exponential decay property $\|u-\bar u_0\|_{H^1}\le Ce^{-a_0 t}\|u_0\|_{H^4}$ as in Theorem 1.3 also holds.
}

\vskip 0.2 true cm

{\bf Remark 1.3. }{\it In Theorem 1.1, the decay property (1.6) comes from the integral
$$
\int_0^t(1+t-\tau)^{-1}(1+\tau)^{-\f{3}{2}}d\tau=O(t^{-1}).
$$
In Theorem 1.2, the decay property (1.10) comes from the integral
$$
\int_0^t(1+t-\tau)^{-\f{3}{4}}(1+\tau)^{-1}d\tau=O(t^{-\f{3}{4}}\ln t)=O(t^{-\f{3}{4}+\dl})\quad\text{for any $\dl>0$}.
$$
One can find more details on the proofs of (1.6) and (1.10) in $\S 5$ and $\S 6$ below, respectively.

In addition, the estimates (1.5) and (1.9) are obtained by the interpolation between the
uniform boundedness of $\|u\|_{C([0, +\infty),H^4)}$ and the decay estimates of the first-order derivatives in (1.4) and (1.8),
respectively.
}

\vskip 0.2 true cm

We now mention some interesting or remarkable results which are closely related to our works. In [6-9], the authors
studied the large time behavior of small perturbed solutions to the compressible isentropic Navier-Stokes equations in the infinite cylindrical domains or parallel domains with Dirichlet boundary conditions of the velocity. For the initial value problem in the whole space, the optimal
decay properties of global small perturbed solutions have been established in [2-6], [11-13] and the references therein.
With respect to the initial-boundary value problem, the optimal decay property have also been obtained in [10], [14-15] and so on.

Next we comment on the proofs of Theorem 1.1-Theorem 1.3.
In [8-9], under the smallness assumptions on the Reynolds number, the Mach
number and the initial perturbation, the authors obtain the large time
behavior of the solution and show that the related profile system is described by the 1-D viscous Burgers equation.
Moreover, it is proved in [8-9] that the decay properties of the solution are determined by the lower-energy part
meanwhile the asymptotic behavior depends on the first simple zero eigenvalue of the resulting linearized operator.
However, in our cases, there are some differences from those in [8-9].
The first difference is: due to the influences of the partial periodic variables in (1.3), we have to use the Fourier series
and Fourier transformation simultaneously to  analyze the linearized equation of (1.3), and which yields that the first zero eigenvalue of
the linearized equation in the low frequency is  of four  multiplicities, while the  zero-eigenvalue is simple and
only Fourier transformation is applied in [8-9].
This leads to that we should consider the difference of the solution $u$ and its mean value  $\bar u$ on
the periodic variables (fortunately, $\bar u$ is just only the right eigenvector of zero eigenvalue) in order to obtain more precise decay property.
Here we point out that only the difference
of the density $\rho$ and its mean value $\bar\rho$ on the cross section needs to be considered in [8-9].
The second difference is: the resulting profile systems (1.7) and (1.11) are the nonlinear partial differential systems
other than the scalar viscous Burgers equation as in [8-9]. Thanks to some delicate analysis on
the linearized system together with a variant
of the Matrumura-Nishida energy method, we can complete the proofs of Theorem 1.1-Theorem 1.3.

\vskip 0.2 true cm

The rest of the paper is organized as follows: In $\S 2$, we will establish some large time properties on
the solution to the homogeneous linearized  problem of (1.3).
In $\S 3$, we study the well-posedness and the asymptotic behavior of the solution to the auxiliary problem (1.7).
In $\S 4$ and $\S 5$, we complete the proof of Theorem 1.1. Finally, Theorem 1.2 and Theorem 1.3
are proved in $\S 6$.

\vskip 0.4 true cm \centerline{\bf \S2. Some analysis on the homogeneous linearized problem of (1.3)}
\vskip 0.3 true cm
In this section, we will study the properties of the solution to the homogeneous linearized  problem of (1.3).
Namely, we consider the following linear problem
\begin{equation}
\left\{
\begin{aligned}
&\p_t v+Lv=0,\\
&v(0,z)=u_0(z),
\end{aligned}
\right.\tag{2.1}
\end{equation}
where $v=(\phi,m)$, $u_0=(\phi_0,m_0)$, and
$
L= \left(
\begin{array}{cc}
0 & \gamma \Div \\
\g\na & -\nu_1\Delta I_3 -\nu_2\na\Div
\end{array}
\right)
$
with $I_3$ being the $3\times3$ identity matrix.

Taking both the Fourier transformation on the $y$ variables and the Fourier series expansion on the $x$ variables in (2.1), we have
\begin{equation}
\left\{
\begin{aligned}
&\ds\f{d}{dt}\left(
\begin{array}{cccc}\hat\phi(t,k,\xi)\\ \hat m(t,k,\xi)\end{array}
\right)+\hat L(k,\xi)\left(
\begin{array}{cccc}\hat\phi(t,k,\xi)\\ \hat m(t,k,\xi)\end{array}
\right)=0,\\
&(\hat\phi, \hat m)(0,k,\xi)=(\hat\phi_0, \hat m_0)(k,\xi),
\end{aligned}
\right.\tag{2.2}
\end{equation}
where $k\in\Bbb Z^{\ell}$, $\xi\in\Bbb R^{3-\ell}$, and
$$
\hat L(k,\xi)=\left(
\begin{array}{cccc} 0&i\g k &i\g\xi\\i\g k^T&\nu_1(|k|^2+\xi^2)+\nu_2k^T k&\nu_2k^T\xi \\i\g\xi^T &\nu_2\xi k^T&(\nu_1|k|^2+(\nu_1+\nu_2)|\xi|^2)I_2\end{array}
\right). \eqno{(2.3)}
$$

It follows from a direct computation that the eigenvalues of $\hat L(k,\xi)$ are
$$
\la_1=\la_2=\nu_1p,\ \la_\pm=\f{1}{2}\{(\nu_1+\nu_2)p\pm\sqrt{(\nu_1+\nu_2)^2p^2-4\g^2p}\}, \eqno{(2.4)}
$$
where $p=|k|^2+|\xi|^2$. Obviously, $\la_1=\la_2=\la_\pm=0$ for $p=0$, which means that
$0$ is the $4$ root of $\hat L(k,\xi)$ as $p=0$. In addition,
$\la_1$ and $\la_2$ are positive numbers for $p>0$,  $\la_\pm$ are complex numbers
when  $p>0$ is small, and $\la_+\rightarrow+\infty$, $\la_-\rightarrow\ds\f{\g^2}{\nu_1+\nu_2}$
as $p\to +\infty$. Therefore, if we choose a fixed positive number $r_0$ with $0< r_0^2 < min\{1,\ds\f{\gamma^2}{(\nu_1+\nu_2)^2}\}$,
then we can conclude that for $p\ge r_0^2$,
$$
\max\{\lambda_1,\lambda_2,\re\lambda_-,\re\lambda_+\}\geq\min\{\nu_1r_0^2,\f{\gamma^2}{2(\nu_1+\nu_2)}\}\triangleq a>0.\eqno{(2.5)}
$$

From this, one can check that the solution to (2.1) can be written as
$$
v(t,z)=\mathcal {U}(t)u_0\triangleq(e^{-tL}u_0)(t,z)=\f{1}{(2\pi)^\ell}\ds\sum_{k\in\Bbb Z^l}e^{ik\cdot x}\int_{\Bbb R^{\Bbb N-\ell}}e^{i\xi\cdot y}e^{-t\hat L(k,\xi)}\hat u_0(k,\xi)d\xi.
$$

If we denote
$$
\bar v(t,y)=\f{1}{(2\pi)^{3-\ell}}\int_{\Bbb R^{\Bbb N-\ell}}e^{i\xi\cdot y}e^{-t\hat L(0,\xi)}\hat u_0(0,\xi)d\xi,
$$
then it is easy to see that $\bar v(t,y)=(\bar\phi, {\bar m})(t,y)=\bar{\mathcal {U}}(t)\bar u_0$ is a solution to the following linear system
\begin{equation}
\left\{
\begin{aligned}
&\partial_{t}\bar\phi+\g\Div'\bar m'=0,\\
&\p_t\bar m_j-\nu_1\Delta'\bar m_j=0,\ j=1,\cdots,l,\\
&\p_t\bar m'-\nu_1\Delta'\bar m'-\nu_2\nabla'\Div'\bar m'+\g\nabla'\bar\phi=0,\\
&\bar v(0,y)=(\bar\phi_0,\bar m_0)(y),
\end{aligned}
\right.\tag{2.6}
\end{equation}
where $\bar m'=(\bar m_{\ell+1}, \cdots, \bar m_3)$.

In particular, if $\ell=3$, we can obtain from (2.6) that
\begin{equation}
\left\{
\begin{aligned}
&\partial_{t}\bar\phi=0,\\
&\p_t\bar m=0,\\
&(\bar\phi,\bar m)(0,y)=(\bar\phi_0,\bar m_0),
\end{aligned}
\right.\tag{2.7}
\end{equation}
which implies $\bar \phi=\bar\phi_0$ and $\bar m=\bar m_0$. This, together with (2.5), yields
$\|v-\bar u_0\|_{H^s}\le C e^{-a_0t}\|u_0\|_{H^s}$ for $s\ge 0$.

Next we study the decay properties of $v$ in the cases of  $\ell=1$ and $\ell=2$.
Let $r_0$ be a positive constant defined in (2.5), and $\chi(z)$ be a function such that $\hat\chi(k,\xi)=\hat\chi(q)=1_{|q|\leq r_0}(k,\xi)$
with $|q|=\sqrt{k^2+|\xi|^2}$. As in [7-9], we set
\begin{align*}
&\mathcal {U}_{(0)}(t)u_0=\mathcal {F}^{-1}(e^{-t\hat L(k,\xi)}\hat\chi(q)\hat u_0),
\ \mathcal {U}_{(\i)}(t)u_0=\mathcal {F}^{-1}(e^{-t\hat L(k,\xi)}(1-\hat\chi(q)) \hat u_0),\\
&\bar{\mathcal {U}}_{(0)}(t)u_0=\mathcal {F}^{-1}(e^{-t\hat L(0,\xi)}\hat\chi(0,\xi) \hat u_0),\ \bar{\mathcal {U}}_{(\i)}(t)u_0
=\mathcal {F}^{-1}(e^{-t\hat L(0,\xi)}(1-\hat\chi(0,\xi)) \hat u_0).\tag{2.8}
\end{align*}
Then we have following results.

{\bf Lemma 2.1.} {\it For $\ell=1,2$ and $u_0\in H^s(\Omega)\cap L^1(\Omega)$
with $\Omega=\Bbb T^\ell\times\Bbb R^{3-\ell}$, $s\in\Bbb N$, then one has

(i) $\|\p^j_y{\mathcal {U}}_{(0)}(t)u_0\|_{L^2(\Omega)}\leq C_j(1+t)^{-\f{3-\ell}{4}-\f{j}{2}}\|u_0\|_{L^1(\Omega)}$, $j\geq0$.

(ii) $(\mathcal {U}_{(0)}(t)u_0)(t,z)=\bar v(t,y)$, and further $\p^{j}_{x}{\mathcal {U}}_{(0)}(t)u_0=0$ for $j\geq1$.

(iii) $\|\mathcal {U}_{(\i)}(t)u_0\|_{H^s(\Omega)}\leq Ce^{-a_0t}\|u_0\|_{H^s(\Omega)}$.

(iv) $\|\bar{\mathcal {U}}_{(\i)}(t)\bar u_0\|_{H^s(\Bbb R^{3-\ell})}\leq Ce^{-a_0t}\|u_0\|_{H^s(\Omega)}$.

Where $a_0$ is a positive constant such that $\ds\f{a}{2}\le a_0\le a$ holds,  and the constant $a$ is defined in (2.5).}

{\bf Proof. } (i) It follows from the Parseval's identity and a direct computation that
\begin{align*}
\|\p^j_y{\mathcal {U}}_{(0)}(t)u_0\|^2_{L^2(\Omega)}&=\int_{|\xi|\leq r_0}|\xi|^{2j}|e^{-t\hat L(0,\xi)}\hat\chi(0,\xi)\hat u_0(0,\xi)|^2d\xi\\
&\leq\|\hat u_0\|^2_{L^\i(\Omega)}\int_{|\xi|\leq r_0}|\xi|^{2j}|e^{-t\hat L(0,\xi)}|^2d\xi\\
&\leq C\|u_0\|^2_{L^1(\Omega)}\int_{|\xi|\leq r_0}|\xi|^{2j}e^{-2a_0t|\xi|^2}d\xi\\
& \leq |Zu|\le\left\{
\begin{aligned}  C t^{-\f{3-\ell}{2}-j}\|u_0\|^2_{L^1(\Omega)},\\
 C \|u_0\|^2_{L^1(\Omega)}. \end{aligned}
\right.
\end{align*}

(ii) follows  the fact of $ \mathcal {U}_{(0)}(t)u_0 =\mathcal {U}_{(0)}(t)\bar u_0$ directly.

By (2.5), we have for $|k|^2+|\xi|^2>r_0$
$$
|e^{-t\hat L(0,\xi)}|\leq e^{-a t}.
$$
Therefore,
\begin{align*}
\|\p_{x}^{j_1}\p_y^{j_2}\mathcal {U}_{(\i)}(t)u_0\|_{L^2(\Omega)}&=
\|\mathcal {F}^{-1}(e^{-t\hat L(k,\xi)}(1-\hat\chi(p))k^{j_1}\xi^{j_2}\hat u_0)\|_{L^2(\Omega)}\\
&\leq Ce^{-a_0t}\|\p_{x}^{j_1}\p_{y}^{j_2}u_0\|_{L^2(\Omega)},
\end{align*}
which derives (iii) and (iv).\qquad \qquad \qquad \qquad $\square$

\vskip 0.4 true cm \centerline{\bf \S3. The global existence of the solution to the auxiliary problem (1.7)}
\vskip 0.2 true cm

In this section, we shall establish the global existence of the classical solution $\eta(t,y)$ to (1.7)
by the energy method.
For notational convenience, we set
$h_1=\left(
\begin{array}{cccc} -\Div'(w'\otimes w_1) \\-\Div'(w'\otimes w')- \alpha\na'(\sigma^2) \end{array}
\right)$.
Then we have

{\bf Lemma 3.1. }{\it For small $\ve>0$, $\eta_0\in H^s(\Bbb R^2)$  and $\|\eta_0\|_{H^s}\le \ve$,
here $s\in\Bbb N$ and $s\ge 2$, then (1.7)
admits a global solution $\eta(t,y)$ such that
$$
\eta \in C([0,\infty),H^s(\Bbb R^2))\cap C^1([0,\infty),H^{s-2}(\Bbb R^2)) \ \ and \ \ \|\eta\|_{H^s}\le C\ve.
$$
}

{\bf Proof. } {Since the system in (1.7) is like the Navier-Stokes equation (1.1) and the local
well-posedness of the initial data problem for the compressible Navier-Stokes system is well-known,
we can derive that (1.7) is locally well-posed, and we omit the proof details here.
In order to show the global well-posedness of (1.7), we only need to establish the related global energy estimates
by the continuity argument.

Set $T^{j,k}=\p_t^j\p_{y}^k$. By $\int T^{j,k}(1.7)\cdot T^{j,k}\eta dy$ and a direct computation, we have
$$
\f{1}{2}\f{d}{dt} \|T^{j,k}\eta\|_2^2 +\nu_1\|T^{j,k}\nabla'w\|_2^2+\nu_2\|\Div'T^{j,k}w' \|_2^2=(T^{j,k}h_1,T^{j,k}w).\eqno{(3.1)}
$$
In addition, it follows from $\int T^{j,k}(1.7)_1\times T^{j+1,k}\sigma dy$
and $\sum\limits _{i=1}^3 \int T^{j,k}(1.7)_{i+1}\times T^{j+1,k}w_i dy$ that
$$
\|T^{j+1,k}\sigma\|_2^2+\gamma(T^{j,k}\Div'w',T^{j+1,k}\sigma)=0\eqno{(3.2)}
$$
and
\begin{align*}
&\|T^{j+1,k}w\|_2^2+\f{1}{2}\f{d}{dt}(\nu_1\|T^{j,k}\nabla 'w \|_2^2+\nu_2\|T^{j,k}\Div'w'\|_2^2) +\gamma(T^{j,k}\nabla'\sigma,T^{j+1,k}w')\\
&=(T^{j,k}h_1,T^{j+1,k}w).\tag{3.3}
\end{align*}
Note that
$$
(T^{j,k}\nabla '\sigma,T^{j+1,k}w')=-\f{d}{dt}(T^{j,k }\sigma,T^{j,k}\Div 'w')+(T^{j+1,k }\sigma,T^{j,k}\Div 'w').
$$
Together with (3.2)-(3.3) and Young's inequality, this yields
$$
\f{1}{2}\|T^{j+1,k}\sigma\|_2^2+\|T^{j+1,k}w\|_2^2+\f{d}{dt}\biggl[\f{\nu_1}{2}\|T^{j,k}\nabla 'w \|_2^2+\f{\nu_2}{2}\|T^{j,k}\Div ' w'\|_2^2-\gamma(T^{j,k }\sigma,T^{j,k}\Div 'w')\biggr]
$$
$$
\le 2\gamma^2 \|T^{j,k}\Div' w'\|_2^2 +(T^{j,k}h_1,T^{j+1,k}w).\eqno{(3.4)}
$$
Computing $(3.2)\times \ds\f{2\gamma^2}{\nu_2} +(3.4)$, we arrive at
\begin{align*}
&\f{1}{2}\|T^{j+1,k}\sigma\|_2^2+\|T^{j+1,k}w\|_2^2+ \f{2\nu_1\gamma^2}{\nu_2} \|T^{j,k+1}w\|_2^2\\
&\quad +\f{d}{dt}\biggl(\f{\nu_1}{2}\|T^{j,k}\nabla ' w\|_2^2 +\f{\nu_2}{2}\|T^{j,k}\Div ' w'\|_2^2-\gamma(T^{j,k }\sigma,T^{j,k}\Div'w') +\f{\gamma^2}{\nu_2}\|T^{j,k}\eta\|_2^2\biggr)\\
&
\leq\f{2\gamma^2}{\nu_2}(T^{j,k}h_1,T^{j,k}w)+(T^{j,k}h_1,T^{j+1,k}w).\tag{3.5}
\end{align*}
Substituting the inequality $\gamma(T^{j,k}\sigma,T^{j,k}\Div ' w' )\le \f{ \nu_2}{2}\|T^{j,k}\Div'w'\|_2^2+\f{\gamma^2}{2\nu_2}\|T^{j,k}\sigma\|_2^2$
into (3.5), we have
\begin{align*}
&\|T^{j+1,k}\eta\|_2^2+ \|T^{j,k+1}w\|_2^2+
\f{d}{dt}\bigl(\|T^{j,k+1}w\|_2^2+\|T^{j,k}\eta\|_2^2\bigr)\\
&\leq C(  T^{j,k}h_1,T^{j,k}w)+C(T^{j,k}h_1,T^{j+1,k}w).\tag{3.6}
\end{align*}

We only prove Lemma 3.1 for the case $s=2$ since the cases of $s\ge 3$ can be treated anagously.
For $j=k=0$, and $j=0,k=1$, we have from (3.6) respectively
$$
\|\p_t\eta\|_2^2+\| \nabla 'w\|_2^2+
\f{d}{dt}\bigl(\|\nabla 'w\|_2^2+\|\eta\|_2^2\bigr)
\le C(h_1, w)+C(h_1,\p_tw)\eqno{(3.7)}
$$
and
$$
\|\p_t\p_y\eta\|_2^2+\|\p_y^2w\|_2^2+
\f{d}{dt}\bigl(\|\p_y^2w\|_2^2+\|\p_y\eta\|_2^2\bigr)
\le C(\p_yh_1, \p_yw )+C(\p_yh_1,\p_t\p_yw).\eqno{(3.8)}
$$

In addition, it follows from (3.1) with $j=0$ and $k=2$ that
$$
\f{d}{dt}\|\p_y^2\eta\|_2^2+\|\p_y^3 w\|_2^2\le C(\p_y^2h_1,\p_y^2 w).\eqno{(3.9)}
$$
Adding (3.7)-(3.9), and using
$
(\p_yh_1,\p_yw)=-(h_1,\p^2_yw)\text{ and }(\p_y^2h_1,\p_y^2 w)=-(\p_y h_1,\p_y^3 w)$,
we get
\begin{align*}
&\|\p_t \eta \|_{H^1}^2+\|\p_y w \|_{H^2}^2+\f{d}{dt}\|\eta\|_{H^2}^2\\
&\le C(h_1, w)+C(h_1,\p_tw)+C(h_1, \p_y^2w)+C(\p_yh_1, \p_y^3w )+C(\p_yh_1,\p_t\p_yw).\tag{3.10}
\end{align*}
On the other hand, by $\nabla' (1.7)_1+ (1.7)_3$, we have
$$
\p_t \nabla '\sigma +\gamma \nabla'\sigma = -\p_t w '+\nu_1 \Delta' w' +(\nu_2-\gamma)\nabla'\Div'w'-\Div'(w'\otimes w')-\alpha\nabla(\sigma^2).\eqno{(3.11)}
$$
Computing $\int(3.11)\cdot\nabla'\sigma dy$ yields
\begin{align*}
&\gamma \|\nabla '\sigma\|_2^2+\f{d}{dt} \|\nabla '\sigma\|_2^2\\
&\le ( -\p_t w '+\nu_1\Delta'w' +(\nu_2-\gamma)\nabla'\Div'w'-\Div'(w'\otimes w')-\alpha\nabla (\sigma^2),\nabla'\sigma))\\
&\le \f{\gamma}{2} \|\nabla ' \sigma\|_2^2+C( \|\p_t w'\|_2^2+\|\p_y^2 w'\|_2^2)+(-\Div'(w'\otimes w')-\alpha\nabla(\sigma^2),\nabla'\sigma)),
\end{align*}
which derives
\begin{align*}
&\f{\gamma}{2}\|\nabla'\sigma\|_2^2+\f{d}{dt}\|\nabla'\sigma\|_2^2\\
&\leq C(\|\p_t w'\|_2^2+\|\p_y^2 w'\|_2^2)+(-\Div'(w'\otimes w')-\alpha\nabla'(\sigma^2),\nabla'\sigma) ).\tag{3.12}
\end{align*}
Computing $\int\p_y (3.11) \cdot \p_y \nabla '\sigma dy $, we have
\begin{align*}
&\gamma \|\p_y \nabla '\sigma\|_2^2+\f{d}{dt} \|\p_y \nabla '\sigma\|_2^2 \\
&\le ( -\p_y\p_t w '+\nu_1 \p_y\Delta' w' +(\nu_2-\gamma)\p_y\nabla '\Div ' w' -\p_y\Div '(w'\otimes w')-\alpha \p_y\nabla (\sigma^2),\p_y\nabla ' \sigma) )\\
&\le \f{\gamma}{2} \|\p_y \nabla ' \sigma\|_2^2+C( \|\p_y\p_t w'\|_2^2+\|\p_y^3 w'\|_2^2)+(-\p_y\Div '(w'\otimes w')-\alpha \p_y\nabla' (\sigma^2),\p_y\nabla ' \sigma) ),
\end{align*}
which concludes
$$
\f{\gamma}{2} \|\p_y^2 \sigma\|_2^2+\f{d}{dt} \|\p_y^2 \sigma\|_2^2
\le C( \|\p_y\p_t w'\|_2^2+\|\p_y^3 w'\|_2^2)+(-\p_y\Div '(w'\otimes w')-\alpha \p_y\nabla' (\sigma^2),\p_y\nabla ' \sigma) ). \eqno{(3.13)}
$$
Combining (3.10) with (3.12)-(3.13) yields
\begin{align*}
&\|\p_t \eta\|_{H^1}^2+\|\p_y w\|_{H^2}^2+\|\p_y\sigma\|_{H^1}^2+\f{d}{dt}\|\eta\|_{H^2}^2\\
&\le  C(h_1, w)+(h_1,\p_tw)+C(h_1, \p_y^2w)+C(\p_yh_1, \p_y^3w)
+(\p_yh_1,\p_t\p_yw)\\
&\quad +(-\Div '(w'\otimes w')-\alpha \nabla' (\sigma^2),\nabla ' \sigma)+(-\p_y\Div'(w'\otimes w')-\alpha\p_y\nabla'(\sigma^2),\p_y\nabla'\sigma).\tag{3.14}
\end{align*}
Define
$$
N(t)=\|\eta\|_{H^2}^2+ \int_0^t \{\|\p_t \eta\|_{H^1}^2+\|\p_y w\|_{H^2}^2+\|\p_y\sigma\|_{H^1}^2\}d\tau.
$$
Integrating (3.14) with respect to the variable $\tau$ over $(0,t)$, one has
\begin{align*}
N(t)&\le \|\eta_0\|_{H^2}^2 + C\int_0^t \bigl\{(h_1, w)+(h_1,\p_tw)+(h_1, \p_y^2w)+(\p_yh_1, \p_y^3w )+(\p_yh_1,\p_t\p_yw)\\
&+(-\Div'(w'\otimes w')-\alpha\nabla'(\sigma^2),\nabla'\sigma)
+ (-\p_y\Div'(w'\otimes w')-\alpha\p_y\nabla'(\sigma^2),\p_y\nabla'\sigma)\bigr\}d\tau.\tag{3.15}
\end{align*}

Next we treat the terms in the right hand side of (3.15) separately.
For the term $\int_0^t (  h_1, w)d\tau $, by using $(1.7)_1$, we see that
$$
(h_1,w)=(-(\Div'(w'\otimes w_1),w_1)+(-\Div'(w'\otimes w')-\alpha\nabla'(\sigma^2),w').
$$
Noting that
\begin{align*}
(-(\Div'(w'\otimes w_1),w_1)&=\int w_1\cdot(w'\cdot\nabla')w_1dy =\f{1}{2}\int (w'\cdot\nabla')|w_1|^2dy\\
&=-\f{1}{2}\int |w_1|^2\Div'w'dy= \f{1}{2\gamma }\int w_1^2 \p_\tau \sigma dy
\end{align*}
and
\begin{align*}
(-\Div'(w'\otimes w')-\alpha \nabla' (\sigma^2),w')&=-\f{1}{2}\int(\Div'w')|w'|^2dy+\int \alpha \sigma^2\Div ' w' dy  \\
&= \f{1}{2\gamma }\int\p_\tau \sigma |w'|^2 dy -\f{1}{\gamma}\int \alpha \sigma^2 \p_t \sigma dy,
\end{align*}
then we have
$$
\int_0^t (h_1, w)d\tau=J_1+J_2,\eqno{(3.16)}
$$
where
\begin{align*}
&J_1=-\f{\alpha}{\gamma}\int_0^t \int \p_\tau \sigma\sigma^2dyd\tau,\\
&J_2=\f{1}{2\gamma}\int_0^t\int\p_\tau\sigma|w|^2dy d\tau.
\end{align*}
It follows from a direct computation that
\begin{align*}
J_1& \le C (\int \sigma^3(t,y) dy-\int \sigma^3(0,y) dy) \\
&\le C (\|\sigma\|_3^3+\|\sigma_0\|_3^3) \le C\|\sigma_0\|_{H^2}^3+CN^{3/2}(t)\tag{3.17}
\end{align*}
and
$$
J_2=\f{1}{2\gamma} \int_0^t \int \p_\tau (\sigma |w|^2 )dy d\tau -\f{1}{\gamma} \int_0^t \int \sigma w\p_\tau wdy d\tau.\eqno{(3.18)}
$$
For the first term in the right hand side of (3.18), we have
$$
\f{1}{2\gamma} \int_0^t \int \p_\tau (\sigma |w|^2 )dy d\tau \le C \int \sigma_0 |w_0|^2 dy +C   \int \sigma |w|^2 dy  \le C\|\eta_0\|_{H^2}^3+CN^{3/2}(t).\eqno{(3.19)}
$$
For the second term in the right hand side of (3.18), we have by $(1.7)_{2,3,4}$
\begin{align*}
&\f{1}{\gamma} \int_0^t \int \sigma w\p_\tau wdy d\tau\\
& \le C\int_0^t \int\sigma\{\nu_1\Delta'w\cdot w+\nu_2\nabla'\Div'w'\cdot w'-\Div'(w'\otimes w')w'-\alpha\nabla'(\sigma^2)w'\\
&\quad -\Div'(w'\otimes w_1)w_1-\g\nabla'\si\cdot w'\}dyd\tau.\tag{3.20}
\end{align*}
It is noted that
$$
\nu_1\int_0^t \int \sigma   \Delta'w \cdot wdyd\tau =-\nu_1 \int_0^t \int \nabla '(\sigma w )\nabla ' w dyd\tau
 \le C\int_0^t \|\eta\|_{L^\infty} \|\nabla \eta\|_2^2 d\tau \le CN^{\f{3}{2}}(t).
$$
Similarly,
$$
\int_0^t\int\nu_2\sigma\nabla'\Div'w'\cdot w'dyd\tau\le C N^{\f{3}{2}}(t),
$$
$$
\int_0^t \int \sigma div'(w'\otimes w')w 'dyd\tau \le C\int_0^t \|\nabla \eta \|_2\|\eta\|_6^3d\tau \le C\int_0^t \|\nabla \eta \|_2\|\eta\|_2^{2}\|\nabla ' \eta\|_2 d\tau\le CN^2(t).
$$
Then substituting the above three estimates into (3.20) and subsequently combining (3.17)-(3.19) yield
$$
\int_0^t (h_1,w)d\tau \le \|\eta_0\|_{H^2}^2 +CN^\f{3}{2}(t).\eqno{(3.21)}
$$
For the other left terms in the right hand side of (3.15), we have
\begin{align*}
\int_0^t(  h_1,\p_\tau w)d\tau &\le \int_0^t \int |\p_y \eta||\eta|  |\p_\tau w| dy d\tau \\
&\le C \|\eta\|_{L^\infty} \int_0^t(\|\p_y \eta\|_2^2 + \|\p_t w\|_2^2) dy d\tau \\
&\le CN^{3/2}(t)\tag{3.22}
\end{align*}
and
\begin{align*}
\int_0^t(  h_1,\p_y^2 w)d\tau &\le \int_0^t \int |\p_y\eta||\eta||\p_y^2 w|dyd\tau \\
&\le C \|\eta\|_{L^\infty} \int_0^t(\|\p_y\eta\|_2^2 + \|\p_y^2 w\|_2^2) dy d\tau \\
&\le CN^{3/2}(t)\tag{3.23}
\end{align*}
and
\begin{align*}
&\int_0^t(h_1,\p_y \sigma)d\tau \le CN^{3/2}(t),\tag{3.24}\\
&\int_0^t(\p_yh_1,\p_t\p_yw)d\tau\le\int_0^t\int|\p_t\p_y w|(|\eta||\p_y\eta|^2+|\eta|^2|\p^2\eta|)dyd\tau\\
&\qquad \leq C\|\eta\|_{L^\i}\int_0^t\|\p_t\p_y w\|\|\p_y\eta\|_{L^4}^2d\tau+C\|\eta\|_{L^\i}^2\int_0^t\|\p_t\p_y w\|\|\p_y^2\eta\|d\tau\\
&\qquad \le CN^{3/2}(t),\tag{3.25}\\
&\int_0^t(  \p_y h_1,\p_y^3 w )d\tau \le C \int_0^t \int (|\p_y ^2 \eta| |\eta|+|\p_y \eta|^2)  |\p_y^3 w|  dy d\tau \\
&\qquad \le C \|\eta\|_{L^\infty} \int_0^t \int |\p_y ^2 \eta||\p_y^3 w|  dy d\tau+  C \int_0^t \|\p_y^3 w\|_2\|\p_y \eta\|_4^2 d\tau \\
&\qquad \le  C \|\eta\|_{L^\infty} \int_0^t \int |\p_y ^2 \eta||\p_y^3 w|  dy d\tau+  C \int_0^t \|\p_y^3 w\|_2\|\p_y \eta\|_2\|\p_y^2 \eta\|_2 d\tau \\
&\qquad \le CN^{3/2}(t),\tag{3.26}\\
&\int_0^t(  \p_y h_1,\p_y^2 \sigma)d\tau \le C \int_0^t \int (|\p_y ^2 \eta| |\eta|+|\p_y \eta|^2)  |\p_y^2 \sigma|  dy d\tau \\
&\qquad \le C \|\eta\|_{L^\infty} \int_0^t \int |\p_y ^2 \eta||\p_y^2 \sigma|  dy d\tau+  C \int_0^t \|\p_y^2 \sigma\|_2\|\p_y \eta\|_4^2 d\tau \\
&\qquad \le  C \|\eta\|_{L^\infty} \int_0^t \int |\p_y ^2 \eta||\p_y^2 \sigma|  dy d\tau+  C \int_0^t \|\p_y^2 \sigma\|_2\|\p_y \eta\|_2\|\p_y^2 \eta\|_2 d\tau \\
&\qquad \le CN^{3/2}(t).\tag{3.27}
\end{align*}
Therefore, substituting (3.21)-(3.27) into (3.15) yields
$$
N(t)\le C\ve^2 +CN^{3/2}(t),
$$
which implies $N(t)\le C\ve^2$ holds for all $t\in [0, +\infty)$. By the local existence of the solution to
(1.7) and the continuity argument, we know that Lemma 3.1 holds for $s=2$. Analogously, Lemma 3.1 also holds
for $s\ge 3$. \qquad\qquad  \qquad  \qquad  \qquad  \qquad    $\square$

Next we establish the decay property of the solution to (1.7). For notational convenience,
we write (1.7) as
$$
\p_t \eta+A\eta =B(\eta),\eqno{(3.28)}
$$
where
$$
A\eta= \left(
\begin{array}{cccc}  \g \Div 'w' \\
-\nu_1\Delta 'w_1 \\
-\nu_1\Delta 'w'-\nu_2 \nabla'\Div'w'+\g\nabla'\sigma
\end{array}
\right)
$$
and
$$
B(\eta)=
 \left(
\begin{array}{cccc} 0 \\
-\Div'(w'\otimes w_1)\\
-\Div'(w'\otimes w')-\alpha \nabla '(\sigma^2)
\end{array}
\right).
$$
Noting that $\hat {A\eta}=\hat L(0,\xi) \hat \eta$, where $\hat L(0,\xi)$ has been defined in (2.2),
then by Duhamel's principle, we have from (3.28)
$$
\eta(t,y)=e^{-tA}\eta_0 +\int_0^t e^{-(t-\tau)A}B(\eta)(\tau)d\tau.  \eqno{(3.29)}
$$

{\bf Lemma 3.2.} {\it Set $M_0(t)=\sup\limits_{0\le \tau \le t} \{ (1+\tau)^{\f{1}{2}} \|\eta\|_2+(1+\tau)\|\p_y\eta\|_2\}$.
If $\eta_0(y)=(\sigma_0,w_0)(y)\in H^4 (\Bbb R^2)\cap L^1 (\Bbb R^2)$ and $\|\eta_0\|_{H^4\cap L^1}\le\ve$, then we have for small $\ve>0$

(i) $ \|\p_y^j e^{-tA}\eta_0\|_2 \le C_j((1+t)^{-\f{1}{2}-\f{j}{2}}\|\eta_0\|_{L^1}+e^{-a_0t}\|\eta_0\|_{H^j})$,
\quad  $j=0, 1, \cdots, 4$.

(ii) $M_0(t) \le C(\|\eta_0\|_{L^1}+\|\eta_0\|_{H^4})+C(M_0^2(t)+\ve^\f{1}{4}M_0^\f{17}{12}(t))$.

Where $a_0$ is a positive constant such that $\ds\f{a}{2}\le a_0\le a$ holds,  and the constant $a>0$ is defined in (2.5).}

\vskip 0.2 true cm

{\bf Proof.}   (i)  Its proof is completely similar to that for Lemma 2.1 (i) and (iii), we omit it here.

(ii) From (3.29) we see that
$$
\|\p_y^j \eta(t,y)\|_2=\|\p_y^je^{-tA}\eta_0\|_2 +\int_0^t \|e^{-(t-\tau)A}\p_y^jB(\eta)\|_2(\tau)d\tau.\eqno{(3.30)}
$$
Since the first term $\|\p_y^je^{-tA}\eta_0\|_2$ in the right hand side of (3.30) has been treated in (i), we only need to treat the second term.
For the nonlinear term $B(\eta)$, we observe that
$$B(\eta)=\Div'b(\eta),$$
where
$$
b(\eta)=\left(
\begin{array}{cccc} 0 \\ -w'\otimes w_1\\-w'\otimes w'-\alpha \sigma^2 I_2\end{array}
\right).
$$
It follows from (i) that
$$\|e^{-(t-\tau)A}B(\eta)\|_2\le C((1+t-\tau)^{-\f{1}{2}} \|B(\eta)\|_{L^1}+e^{-a_0(t-\tau)}\|B(\eta)\|_{L^2})\eqno{(3.31)}$$
and
$$\|e^{-(t-\tau)A}B(\eta)\|_2 =\|\Div'e^{-(t-\tau)A}b(\eta)\|_2 \le C((1+t-\tau)^{-1} \|b(\eta)\|_{L^1}+e^{-a_0(t-\tau)}\|\p_yb(\eta)\|_{L^2}).\eqno{(3.32)}$$
In addition, we have
$$
\int_0^t \|e^{-(t-\tau)A} B(\eta)\|_2(\tau)d\tau =\int_0^{t/2} \|\p_ye^{-(t-\tau)A} b(\eta)\|_2(\tau)d\tau+\int_{t/2}^t \|e^{-(t-\tau)A} B(\eta)\|_2(\tau)d\tau.\eqno{(3.33)}
$$
For the first term in the right hand side of (3.33), one has by (3.32)
\begin{align*}
&\int_0^{t/2} \|\p_ye^{-(t-\tau)A} b(\eta)\|_2(\tau)d\tau\\
&\le C\int_0^{t/2}\bigl((1+t-\tau)^{-1} \|b(\eta)\|_{L^1}+e^{-a_0(t-\tau)}\|\p_yb(\eta)\|_{L^2}\bigr)d\tau\\
& \le C\int_0^{t/2}\bigl((1+t-\tau)^{-1}\|\eta\|_2^2+e^{-a_0(t-\tau)}\|\eta\|_{L^4}\|\p_y\eta\|_{L^4}\bigr) d\tau\\
&  \le C\int_0^{t/2}\bigl((1+t-\tau)^{-1} (1+\tau)^{-1}M_0^2(\tau)
+\ve^\f{1}{4}e^{-a_0(t-\tau)}(1+\tau)^{-\f{13}{8}}M_0^\f{7}{4}(\tau)\bigr)d\tau\\
& \le C(1+t)^{-1}\ln(1+t)M^2_0(t)+\ve^\f{1}{4}e^{-\f{a_0}{4}t}M_0^\f{7}{4}(t).
\tag{3.34}
\end{align*}
For the second term of (3.33), then it follows from (3.31) that
\begin{align*}
&\int_{t/2}^t \|e^{-(t-\tau)A} B(\eta)\|_2(\tau)d\tau\\
&\le C\int_{t/2}^t\bigl( (1+t-\tau)^{-\f{1}{2}} \|B(\eta)\|_{L^1}+e^{-a_0(t-\tau)}\|B(\eta)\|_{L^2}\bigr) d\tau\\
&\le C\int_{t/2}^t\bigl((1+t-\tau)^{-\f{1}{2}}\|\p_y\eta\|_2\|\|\eta\|_2+e^{-a_0(t-\tau)}\|\eta\|_{L^4}\|\p_y\eta\|_{L^4}\bigr)d\tau\\
&\le C\int_{t/2}^t\bigl((1+t-\tau)^{-\f{1}{2}} (1+\tau)^{-\f32}M_0^2(\tau)
+\ve^\f{1}{4}e^{-a_0(t-\tau)}(1+\tau)^{-\f{13}{8}}M_0^\f{7}{4}(\tau)\bigr) d\tau\\
&\le C (1+t)^{-1}M_0^2(t)+C\ve^\f{1}{4}(1+t)^{-\f{13}{8}}M_0^\f{7}{4}(t).\tag{3.35}
\end{align*}
Therefore, combining (3.34) with (3.35),  we can obtain the estimate (3.30) with $j=0$.

In the case of $j=1$, we have
\begin{align*}
&\int_0^{t/2} \|\p_y^2 e^{-(t-\tau)A} b(\eta)\|_2(\tau)d\tau\\
\le& C\int_0^{t/2}\bigl((1+t-\tau)^{-\f{3}{2}} \|b(\eta)\|_{L^1}+e^{-a_0(t-\tau)}\|\p_y^2b(\eta)\|_{L^2}\bigr) d\tau\\
\le& C\int_0^{t/2}\bigl((1+t-\tau)^{-\f{3}{2}}\|\eta\|_2^2+e^{-a_0(t-\tau)}(\|\p_y\eta\cdot\p_y\eta\|_{L^2}+\|\eta\cdot\p_y^2\eta\|_{L^2})\bigr)d\tau\\
\le& C\int_0^{t/2}\bigl((1+t-\tau)^{-\f{3}{2}} (1+\tau)^{-1}M_0^2(\tau)+e^{-a_0(t-\tau)}(\|\p_y\eta\|^2_{L^4}+\|\eta\|_{L^4}\|\p_y^2\eta\|_{L^4})\bigr)d\tau\\
\le& C\int_0^{t/2}\bigl((1+t-\tau)^{-\f{3}{2}} (1+\tau)^{-1}M_0^2(\tau)+e^{-a_0(t-\tau)}(\ve^\f{1}{2}(1+\tau)^{-\f{3}{2}}M_0^\f{3}{2}(\tau)
+\ve^\f{7}{12}(1+\tau)^{-\f{3}{2}}M_0^\f{17}{12}(\tau))\bigr)d\tau\\
\le& C(1+t)^{-\f{3}{2}}\ln(1+t)M^2_0(t)+\ve^\f{1}{2}e^{-\f{a_0}{4}t}M_0^\f{17}{12}(t) \tag{3.36}
\end{align*}
and
\begin{align*}
&\int_{t/2}^t \| \p_y e^{-(t-\tau)A}   B(\eta)\|_2(\tau)d\tau\\
\le& C\int_{t/2}^t \bigl((1+t-\tau)^{-1} \|B(\eta)\|_{L^1}+e^{-a_0(t-\tau)}\|\p_yB(\eta)\|_{L^2}\bigr)d\tau\\
\le& C\int_{t/2}^t\bigl((1+t-\tau)^{-1}\|\p_y\eta\|_2\|\|\eta\|_2+e^{-a_0(t-\tau)}(\|\p_y\eta\cdot\p_y\eta\|_{L^2}
+\|\eta\cdot\p_y^2\eta\|_{L^2})\bigr)d\tau\\
\le& C\int_{t/2}^t\bigl((1+t-\tau)^{-1} (1+\tau)^{-\f32} M_0^2(\tau)+e^{-a_0(t-\tau)}(\ve^\f{1}{2}(1+\tau)^{-\f{3}{2}}M_0^\f{3}{2}(\tau)
+\ve^\f{7}{12}(1+\tau)^{-\f{3}{2}}M_0^\f{17}{12}(\tau))\bigr)d\tau\\
\le& C ((1+t)^{-\f{3}{2}}\ln(1+t)M^2_0(t)+\ve^\f{1}{2}(1+\tau)^{-\f{3}{2}}M_0^\f{3}{2}(t)
+\ve^\f{7}{12}(1+\tau)^{-\f{3}{2}}M_0^\f{17}{12}(t)).\tag{3.37}
\end{align*}

Combining (3.34)-(3.37) yields for $j=0, 1$,
$$
\|\p_y ^j \eta\|\le C(1+t)^{-\f{1}{2}-\f{j}{2}}\|\eta_0\|_{H^4\cap L^1}
+C(1+t)^{-\f{1}{2}-\f{j}{2}}(M^2_0(t)+\ve^\f{1}{4}M_0^\f{17}{12}(t)).
$$
Consequently, we complete  the proof of (ii).\qquad \qquad \qquad $\square$

\vskip 0.2 true cm

{\bf Remark 3.1.} {\it By analogous arguments as in Lemma 3.2, for the problem (1.11), if
$\eta_0(y)=(\sigma_0,w_0)(y)\in H^4 (\Bbb R)\cap L^1 (\Bbb R)$ and $\|\eta_0\|_{H^4\cap L^1}\le\ve$,
then we have for small $\ve>0$
$${\t M}_0(t) \le C(\|\eta_0\|_{L^1}+\|\eta_0\|_{H^4})+C({\t M}_0^2(t)+\ve^\f{1}{4}{\t M}_0^\f{17}{12}(t)),\eqno{(3.38)}$$
where ${\t M}_0(t)=\sup\limits_{0\le \tau \le t} \{ (1+\tau)^{\f{1}{4}} \|\eta\|_2+(1+\tau)^{\f{3}{4}}\|\p_y\eta\|_2\}$.}

\vskip 0.4 true cm \centerline{\bf \S4. Some decay properties of the solution $u(t,z)$ to (1.3) for $z\in\Bbb T \times \Bbb R^2$}
\vskip 0.2 true cm

In this section, we will study the large time behavior of the solution to (1.3) for $z\in\Bbb T \times \Bbb R^2$
as $t$ tends to infinity.
As an ingredient of the proof, we require to cite the following Gagliardo-Nirenberg-Sobolev inequalities
(see [1] and so on) which will be used repeatedly.

{\bf Lemma 4.1.(Gagliardo-Nirenberg-Sobolev)} { \it (i) Let $2\leq p\leq\infty$, and let $j$ and $k$ be the integers satisfying
$$
0\leq j<k,\  k > j+n(\f{1}{2}-\f{1}{p}).
$$
Then for $f\in H^k(\Bbb R^n)$, there exists a constant $C>0$ independent of $f$ such  that
$$
\|\p_x^j f\|_{L^p(\Bbb R^n)}\leq C\|f\|_{L^2(\Bbb R^n)}^{1-a}\|\p_x^kf\|_{L^2(\Bbb R^n)}^a,
$$
where $a=\ds\f{1}{k}(j+\f{n}{2}-\f{n}{p})$.

(ii) If $f\in H^1(\Bbb R)$, then
$$
\|f\|_{L^\i}\leq C\|f\|_{L^2}^{\f{1}{2}}\|\p_x f\|_{L^2}^\f{1}{2}.
$$}

Define
$$
\Pi [g](y)=\bar g(y) =\f{1}{2\pi}\int_{\Bbb T}g(x,y)dx.
$$
By Duhamel's principle, one can get the follow expressions of the solutions to (1.3)
$$
u(t,z)=\mathcal {U}(t)u_0+\int_0^t\mathcal {U}(t-\tau)f(\tau,z)d\tau,\eqno{(4.1)}
$$
and
$$
\bar u(t,y)=\bar{\mathcal {U}}(t)\bar u_0+\int_0^t\bar{\mathcal {U}}(t-\tau)\bar f(\tau,y)d\tau,\eqno{(4.2)}
$$
where $f(t,z)=\left(
\begin{array}{cccc} 0\\G(\phi,m)\end{array}
\right)$ and $\bar f(t,y)=\Pi [f](t,y)$.

We note that $G(\phi,m)$ in (1.3) can be written as
$$
G=\Div{\mathcal {G}}(\phi,m)+\na g(\phi,m),
$$
where
\begin{align*}
\mathcal {G}(\phi,m)&=-\f{\g}{\phi+\g}m\otimes m-\nu_1\na\otimes(\f{\phi}{\phi+\g}m),\\
g(\phi,m)&=-\nu_2\Div(\f{\phi}{\phi+\g}m)-F(\phi),
\end{align*}
and $F(\phi)=\ds\f{\phi^2}{\g^2}\int_0^1(1-\th)^2p^{''}(1+\f{\th\phi}{\g})d\th$.

Due to
$
\Pi[\Div v]=\Pi[ \Div ' v']$ and $\Pi[\na f]=(0,\nabla'\bar f)^T$,
$\bar u(t,y)=(\bar\phi,\bar m)(t,y)$ satisfies the nonlinear system
\begin{equation}
\left\{
\begin{aligned}
&\partial_{t}\bar\phi+\g \Div '\bar m'=0,\\
&\p_t\bar m_1-\nu_1\Delta '\bar m_1=- \Pi[\Div' (\f{\gamma}{\phi+\gamma}m_1 m')-\nu_1 \Delta ' (\f{\phi}{\phi+\gamma}m_1)]\\
&\p_t\bar m'-\nu_1\Delta '\bar m'-\nu_2 \nabla'\Div'\bar m'+\g\nabla'\bar\phi= \Pi[G'(\phi,m')],\\
&\bar u(0,y)=(\bar\phi_0,\bar m_0)(y),
\end{aligned}
\right.\tag{4.3}
\end{equation}
where
$$
G'(\phi,m')=-\Div'(\f{\g}{\phi+\g}m'\otimes m')-\nu_1\Dl'(\f{\phi}{\phi+\g}m')-\nu_2\na'\Div'(\f{\phi}{\phi+\g}m')-\na'F(\phi).
$$

Let $\tilde u=u-\bar u$ and
$$
M(t)=M[u](t)\equiv\sup_{0\leq\tau\leq t}\{(1+\tau)^{\f{1}{2}}\|u(\tau,z)\|_{L^2(\Omega)}+(1+\tau)\|\na u(\tau,z)\|_{L^2(\Omega)}\}.
$$
Noting that $\bar u(t,y)$ is the mean value of $u(t,z)$ with respect to the periodic variable $x$, then one has
$M[\bar u](t)\le CM(t)$ and further $M[\t u](t)\le CM(t)$.

We now establish the following decay estimates on the solution $u$ to (1.3).

{\bf Lemma 4.2. }{\it Assume $u\in C([0, +\infty), H^4(\Bbb T\times\Bbb R^2)\cap L^1(\Bbb T\times\Bbb R^2))$,
$\|u\|_{H^4\cap L^1}\leq\ve$ and $M(t) \le 1$, then for small $\ve>0$

(i) $\|u\|_{L^\i}\le C\ve^\f{1}{6}(1+t)^{-\f{7}{12}}M^\f{5}{8}(t)$.

(ii) $\|  \mathcal {G}(\phi,m)\|_{L^1}\leq C(1+t)^{-1}M^2(t)$.

(iii) $\| g(\phi,m)\|_{L^1}\leq C(1+t)^{-1}M^2(t)$.

(iv) $\|G\|_{H^1} \le C\ve^\f{13}{15}  (1+t)^{-1}M(t)+C\ve^\f{1}{6}(1+t)^{-\f{13}{12}}M^\f{5}{4}(t)$.

(v) $\|G\|_{L^1} \le C\ve^\f{1}{3}(1+t)^{-\f{7}{6}}M^{\f{4}{3}}(t)$.
}

{\bf Proof. }  It follows from Lemma 4.1 and a direct computation that
\begin{align*}
\|u\|_{L^\i}&\le \|\bar u\|_{L^\i}+\|\tilde u\|_{L^\i}\\
& \le C \|\bar u \|_2^\f{1}{2}\|D^2 \bar u \|_2^\f{1}{2} +C\|\tilde u \|_2^\f{5}{8}\|D^4 \tilde u \|_2^\f{3}{8}\\
& \le C \|\bar u \|_2^\f{1}{2} \|D\tilde u \|_2^\f{1}{3}\|D^4 \bar u\|_2^\f{1}{6}+C\|\tilde u \|_2^\f{5}{8}\|D^4 \tilde u \|_2^\f{3}{8}\\
&\le C (1+t)^{-\f{7}{12}}M^\f{5}{6}(t)\ve^\f{1}{6} +C(1+t)^{-\f{5}{8}}M^\f{5}{8}(t)\ve^\f{3}{8}\\
& \le C \ve^\f{1}{6}(1+t)^{-\f{7}{12}}M^\f{5}{8}(t),
\end{align*}
which means that (i) holds.

Since
$$
\|u\cdot u\|_{L^1}\leq\| u\|_{L^2}\|u\|_{L^2}\le C(1+t)^{-1}M^2(t)
$$
and
$$
\|u\cdot D u\|_{L^1}\leq \| u\|_{L^2}\|Du\|_{L^2}\le C(1+t)^{-\f{3}{2}}M^2(t),
$$
we arrive at
$$
\|\mathcal {G}(\phi,m)\|_{L^1} \le C\|u\|_{H^1}^2 \|u\|_{L^2}\le C(1+t)^{-1}M^2(t)
$$
and
$$
\| g (\phi,m)\|_{L^1} \le C\|u\|_{H^1}^2 \|u\|_{L^2}\le C(1+t)^{-1}M^2(t).
$$
Thus, (ii) and (iii) are proved.

Next we show (iv).  Note that
\begin{align*}
&\|u\cdot u\|_{L^2}\leq C\|u\|_{L^4}^2
\le C\|u\|^\f{1}{4}_{L^2}\|Du\|^\f{3}{4}_{L^2}\|u\|^\f{1}{4}_{L^2}\|Du\|^\f{3}{4}_{L^2}\\
&\qquad\qquad  \le C(1+t)^{-\f{7}{4}}M^2(t),\tag{4.4}\\
&\|D(u\cdot u)\|_{L^2}\leq C\|Du\cdot u\|_{L^2}\leq C\|Du\|_{L^2}\|u\|_{L^\i}\\
&\qquad\qquad  \le C\ve^\f{1}{6}(1+t)^{-\f{19}{12}}M^\f{13}{8}(t),\tag{4.5}\\
&\|D^2(u\cdot u)\|_{L^2}\le C\|D^2u\cdot u\|_{L^2}+C\|Du\cdot D u\|_{L^2}\le C\|D^2u\|_{L^2}\|u\|_{L^\i}+C\|D u\|_{L^4}^2\\
&\qquad\qquad  \le C\|u\|_{L^\i}\|Du\|_{L^2}^\f{1}{2}\|D^3u\|_{L^2}^\f{1}{2}+C\|D u\|_{L^2}^{\f{5}{4} } \|D^3u\|_{L^2}^\f{3}{4}\\
&\qquad\qquad  \le C\ve^\f{2}{3}(1+t)^{-\f{13}{12}}M^\f{5}{4}(t)\tag{4.6}
\end{align*}
and
\begin{align*}
\|D^3(u\cdot u)\|_{L^2}&\le C\|u\|_{L^5}\|D^3u\|_{L^\f{10}{3}}+C\|D^2u\|_{L^2}\|D u\|_{L^\i}\\
&\le C\|u\|_{L^2}^\f{1}{10}\|Du\|_{L^2}^\f{9}{10}\|Du\|_{L^2}^\f{2}{15}\|D^4u\|_{L^2}^\f{13}{15}+C\|D u\|_{L^2}^\f{1}{2}\|D^3u\|_{L^2}^\f{1}{2}\|Du\|_{L^2}^\f{1}{2}\|D^4u\|_{L^2}^\f{1}{2}\\
&\le C\ve^\f{13}{15}(1+t)^{-\f{13}{12}}M^\f{17}{15}(t)+C\ve(1+t)^{-1}M(t)\\
&\le C\ve^\f{13}{15}(1+t)^{-1}M(t).\tag{4.7}
\end{align*}

Thus, by (4.4)-(4.6) we have
$$
\|\Div(\f{\g}{\phi+\g}m\otimes m)+\na F(\phi) \|_{L^2}\le C\|D(u\cdot u)\|_{L^2}+C\|D\phi\cdot u\cdot u\|_{L^2}\le C\ve^\f{1}{6}(1+t)^{-\f{19}{12}}M^\f{13}{8}(t)\eqno{(4.8)}
$$
and
$$
\|\nu_1\Dl(\f{\phi}{\phi+\g}m)+\nu_2\na div(\f{\phi}{\phi+\g}m)\|_{L^2}\le C\ve^\f{2}{3}(1+t)^{-\f{13}{12}}M^\f{5}{4}(t).\eqno{(4.9)}
$$

Since
$$
G(\phi,m)=-\Div(\f{\g}{\phi+\g}m\otimes m)-\na F(\phi)-\nu_1\Dl(\f{\phi}{\phi+\g}m)-\nu_2\na\Div(\f{\phi}{\phi+\g}m),\eqno{(4.10)}
$$
then substituting (4.8)-(4.9) into (4.10) yields
$$
\|G\|_{L^2}\le C\ve^\f{1}{6}(1+t)^{-\f{13}{12}}M^\f{5}{4}(t).\eqno{(4.11)}
$$
Similarly,  applying (4.7) and (4.11), we have
$$
\|G\|_{H^1} \le C\|D^3(u\cdot u)\|_{L^2} + C\|G\|_{L^2} \le C\ve^\f{13}{15} (1+t)^{-1}M(t) + C\ve^\f{1}{6}(1+t)^{-\f{13}{12}}M^\f{5}{4}(t).
$$

Finally, we prove (v). It follows a direct computation that
\begin{align*}
\|G\|_{L^1} &\le C \|D(u\cdot u)\|_{L^1} +C\|D^2 (u\cdot u)\|_{L^1} \le C\|u\|_{L^2}\|Du\|_{L^2}+C\|Du\|_{L^2}^2+C\|D^2 u\|_{L^2}\|u\|_{L^2}\\
& \le C(1+t)^{-\f{3}{2}}M^2(t)+C(1+t)^{-2}M^2(t)+C\|Du\|_{L^2}^{\f{2}{3}}\|D^4 u\|_{L^2}^\f{1}{3}\|u\|_{L^2}\\
& \le C (1+t)^{-\f{3}{2}}M^2(t)+C(1+t)^{-2}M^2(t)+C\ve^\f{1}{3}(1+t)^{-\f{7}{6}}M^{\f{4}{3}}(t)\\
&\le C \ve^\f{1}{3}(1+t)^{-\f{7}{6}}M^{\f{4}{3}}(t).
\end{align*}
\qquad \qquad \qquad \qquad \qquad \qquad \qquad \qquad \qquad
\qquad \qquad \qquad \qquad $\square$

For later uses, we specially list some intermediate results in the proof procedure of Lemma 4.2 as follows:

{\bf Lemma 4.3. }{\it  Assuming $u\in C([0, +\infty), H^4(\Bbb T\times\Bbb R^2)\cap L^1(\Bbb T\times\Bbb R^2))$, $\|u\|_{H^4\cap L^1}\leq\ve$ and $M(t) \le 1$, then we have

(i) $\|u\cdot u\|_{L^2}\le C(1+t)^{-\f{7}{4}}M^2(t)$.

(ii) $\|D (u \cdot u) \|_{L^2} \le C\ve^{\f{1}{6}} (1+t)^{-\f{19}{12}}M^{\f{13}{8}}(t)$.

(iii) $\|D^2(u\cdot u)\|_{L^2}\le C\ve^\f{2}{3}(1+t)^{-\f{13}{12}}M^\f{5}{4}(t)$.

(iv) $ \|D^3(u\cdot u)\|_{L^2}\le C\ve^\f{13}{15}(1+t)^{-1}M(t)$.
}

\vskip 0.2 true cm
Next we establish some estimates on $M[\bar u](t)$.

{\bf Lemma 4.4. }{\it Assuming $u\in C([0, +\infty), H^4(\Bbb T\times\Bbb R^2)\cap L^1(\Bbb T\times\Bbb R^2))$,
$\|u\|_{H^4\cap L^1}\leq\ve$ and $M(t) \le 1$, then one has
$$
M[\bar u](t) \le  C\ve +C\ve^\f{1}{6}  M(t)+ CM^{\f{4}{3}}(t).\eqno{(4.12)}
$$
}

{\bf Proof. } By (4.2) and Lemma 2.1, in order to show (4.12), we only need to estimate the integral
\begin{align*}
\int_0^t\bar{\mathcal {U}}(t-\tau)\bar f(\tau,y)d\tau&=\int_0^t\bar{\mathcal {U}}_{(0)}(t-\tau)\bar f(\tau,y)d\tau
+\int_0^t\bar{\mathcal {U}}_{(\i)}(t-\tau)\bar f(\tau,y)d\tau\\
&\triangleq J_1(t,y)+J_2(t,y).
\end{align*}

For $J_1(t,y)$, by Lemma 2.1 (ii), Lemma 4.2 (ii)-(iii) and the expression of $G$, we have for $j=0,1$,
\begin{align*}
\|\p_y^j&J_1(t,y)\|_{L^2}\le C\int_{\f{t}{2}}^t\|\mathcal {U}_{(0)}(t-\tau)\p_y^{j+1} ( \bar {\mathcal {G}}+\bar g)\|_{L^2}d\tau +
C\int_0^{\f{t}{2}}\|\p_y^{j}\mathcal {U}_{(0)}(t-\tau)  \bar G \|_{L^2}d\tau\\
&\le C\int_{\f{t}{2}}^t(1+t-\tau)^{-\f{j+2}{2}}\|\bar {\mathcal {G}}+\bar g\|_{L^1}d\tau +C\int_0^{\f{t}{2}}(1+t-\tau)^{-\f{j+1}{2}}\|\bar G\|_{L^1}d\tau\\
&\le C (1+t)^{-1}M^2(t) \int_{\f{t}{2}}^t (1+t-\tau)^{-\f{j+2}{2}}d\tau +C\ve^{\f{1}{3}} (1+t)^{-\f{j+1}{2}} M^{\f{4}{3}}(t)\int_0^{\f{t}{2}} (1+\tau)^{-\f{7}{6}}d\tau\\
&\le C(1+t)^{-\f{j+1}{2}}M^{\f{4}{3}}(t).
\end{align*}
On the other hand, we have
$$
\|J_2(t,y)\|_{H^1}\le CM(t)\int_0^te^{-a_0(t-\tau)} \ve^\f{1}{6} (1+\tau)^{-1} d\tau
\le C\ve^\f{1}{6} (1+t)^{-1}M (t).
$$
Consequently,
$$
\|\bar u\|_2\le C(1+t)^{-\f{1}{2}}\|\bar u_0\|_2  +C\ve^\f{1}{6} (1+t)^{-1}M(t)+(1+t)^{-\f{1}{2}}M^{\f{4}{3}}(t)
$$
and
$$
\|D\bar u\|_2\le C(1+t)^{-1}\|D\bar u_0\|_2 +C\ve^\f{1}{6} (1+t)^{-1}M(t)+(1+t)^{-1}M^{\f{4}{3}}(t).
$$
Then
$$
M[\bar u ](t)\le C\ve +C\ve^\f{1}{6}  M(t)+CM^{\f{4}{3}}(t)
$$
and  the proof of Lemma 4.4 is completed.\qquad \qquad \qquad \qquad \qquad $\square$

Next, we treat the difference $\tilde u=u-\bar u$. It is noted that $\tilde u$ satisfies the following system
\begin{equation*}
\left\{
\begin{aligned}
&\p_t \tilde u +L\tilde u=f(\phi,m)-\bar f(\phi,m)\triangleq\tilde f(t,z),\\
&\tilde u(0,z)=\tilde u_0(z)=u_0(z)-\bar u_0(y).
\end{aligned}
\right.
\end{equation*}
then by Duhamel's principle, we have
$$
\tilde u(t,z)=\mathcal {U}(t)\tilde u_0 +\int_0^t \mathcal {U}(t-\tau) \tilde f(\tau,z) d\tau .
$$
Using Lemma 2.1, we see that
$$
\|\tilde u\|_{H^1} \le Ce^{-a_0t}\|\tilde u_0\|_{H^1}+ C\int_0^t e^{a_0(t-\tau)}\|\tilde f(\tau,z)\|_{H^1}d\tau.
$$
Furthermore, we have

{\bf Lemma 4.5. }{\it Assuming  $u\in C([0, +\infty), H^4(\Bbb T\times\Bbb R^2)\cap L^1(\Bbb T\times\Bbb R^2))$,
$\|u\|_{H^4\cap L^1}\leq\ve$ and $ M(t) \le 1$, then
$$
(1+t)\|\tilde u\|_{H^1} \le C(\ve+\ve^\f{13}{15}M(t)+\ve^\f{1}{6}M^\f{5}{4}(t)).
$$
}

{\bf Proof. }{By Lemma 4.2, we have
\begin{align*}
\|\tilde u\|_{H^1} &\le Ce^{-a_0t}\|\tilde u_0\|_{H^1}+ C\int_0^t e^{-a_0(t-\tau)}\|\tilde f(\tau,z)\|_{H^1}d\tau\\
&\le Ce^{-a_0t}\|\tilde u_0\|_{H^1}+C\int_0^t e^{-a_0(t-\tau)}\{\ve^\f{13}{15}(1+\tau)^{-1}M(\tau)+C\ve^\f{1}{6}(1+t)^{-\f{13}{12}}M^\f{5}{4}(\tau)\}d\tau\\
& \le Ce^{-a_0t}\|\tilde u_0\|_{H^1}+C\ve^\f{13}{15}(1+t)^{-1}M(t) +C\ve^\f{1}{6}(1+t)^{-\f{13}{12}}M^\f{5}{4}(t)
\end{align*}
and further
$$
(1+t)\|\tilde u\|_{H^1} \le C(\ve+\ve^\f{13}{15} M(t)+\ve^\f{1}{6}M^\f{5}{4}(t)).
$$
This completes the proof of Lemma 4.5.\qquad \qquad \qquad $\square$
}

Finally, we conclude the decay results of $u(t,z)$ in the section.

{\bf Proposition 4.6. }{\it For small $\ve>0$, if $u_0\in H^4(\Bbb T \times \Bbb R^2)$
and $\|u_0\|_{H^4\cap L^1}\leq\ve$, then we have
$$ \|\p_z^k u(t,z)\|_{L^2} = O(t^{-\f{1}{2}-\f{k}{2}}),\ \ k=0,1,\eqno{(4.13)}$$
$$ \|\p_z^k( u(t,z)-\bar u(t,y))\|_{L^2} = O(t^{-1}),\ \ k=0,1,\eqno{(4.14)}$$
where $u$ is the solution of (1.3) and $\bar u$ is the solution of the system (4.3).
}

{\bf Proof. }{From Lemma 4.4-4.5, we see that
$$
M(t)\le C(\ve+\ve^\f{1}{6} M(t)+\ve^\f{1}{6}M^\f{5}{4}(t)),
$$
which implies $M(t) \le C\ve$ for small $\ve>0$ and further
$$\|\p_z^k u(t,z)\|_{L^2} = O(t^{-\f{1}{2}-\f{k}{2}}),\ \ k=0,1.$$
In addition, by Lemma 4.5, we can easily obtain
$$ \|\p_z^k( u(t,z)-\bar u(t,y))\|_{L^2} = O(t^{-1}),\ \ k=0,1.$$
Thus Proposition 4.6 is proved.\qquad\qquad\qquad $\square$
}

\vskip 0.4 true cm \centerline{\bf \S5. Large-time behavior of the solution to (1.7) and the proof of Theorem 1.1}
\vskip 0.3 true cm
In this section, we will mainly complete the proof of Theorem 1.1. To this end,
at first we study the asymptotic behavior of the solution to (1.7).
Let $\eta$ and $\bar u $ be the solutions of (1.7) and (4.3), respectively. We define
$$N_1(t)=\sup \limits_{0\le \tau\le t}(1+\tau) \|\bar u(t,y) -\eta(t,y)\|_{H^1(\Bbb R^2)}$$
and
$$
h=\left(
\begin{array}{cccc} 0 \\h_1\end{array}
\right)=\left(
\begin{array}{cccc}  0\\ -\Div' ( w_1 w') \\-\Div'(w'\otimes w')- \alpha\na'(\sigma^2) \end{array}
\right).
$$

{\bf Lemma 5.1. }{\it   If $u_0\in H^4(\Bbb T \times \Bbb R^2)$ and $\|u_0\|_{H^4\cap L^1}\le \ve$, then for small $\ve>0$
$$
\|\bar u -\eta \|_{H^1} \le C (1+t)^{-1} \ve.
$$
}

{\bf Proof. }{By Duhamel's principle, one has
\begin{align*}
(\bar u-\eta)(t,y)&=\int_0^t\bar{\mathcal {U}}(t-\tau)(\bar f-h)(\tau,y) d\tau \\
&=\int_0^t\bar{\mathcal {U}}_{(0)}(t-\tau)(\bar f-h)  d\tau+\int_0^t\bar{\mathcal {U}}_{(\infty)}(t-\tau)(\bar f-h)  d\tau\\
&=J_1(t,y)+J_2(t,y).\tag{5.1}
\end{align*}
By Proposition 4.6 and Lemma 3.2, we see that
$M(t)+M_0(t)\le C \ve,$
where $M(t)$ and $M_0(t)$
are defined in $\S 4$ and Lemma 3.2 respectively. This, together with Lemma 4.2 (iv), yields
$$
\|\bar f-h\|_{H^1}\le C\|G\|_{H^1}+C\|h_1\|_{H^1}\le C \ve^2(1+t)^{-1}.\eqno{(5.2)}
$$

We now treat the terms $J_1(t,y)$ and $J_2(t,y)$ in the right hand side of (5.1) separately.
At first, we deal with $J_2(t,y)$. It follows from a direct computation that
\begin{align*}
\|J_2\|_{H^1} &\le C\int_0^t e^{-a_0(t-\tau)} \| \bar f-h\|_{H^1}d\tau \\
& \le C\int_0^t e^{-a_0(t-\tau)}\ve^2(1+\tau)^{-1}d\tau \\
&\le C\ve^2(1+t)^{-1}.\tag{5.3}
\end{align*}
Next, we deal with $J_1(t,y)$. It is noted that $\bar f$ can be written as
\begin{align*}
&\bar f =\left(
\begin{array}{cccc} 0 \\ \Pi[-\Div' (\f{\gamma}{\gamma+\phi}m_1m')]\\ \Pi[-\Div' (\f{\gamma}{\gamma+\phi}m'\otimes m')-\nabla ' F(\phi)]
\end{array}
\right)
+ \left(
\begin{array}{cccc} 0 \\ \Pi[-\nu_1\Delta' (\f{\phi}{\gamma+\phi}m_1)]\\ \Pi[-\nu_1\Delta'(\f{\phi}{\gamma+\phi}m')-\nu_2\nabla '\Div '(\f{\phi}{\gamma+\phi}m')]\end{array}
\right)\\
&=\bar f_1+\bar f_2.\\
\end{align*}
From Lemma 2.1 (i), we see that
\begin{align*}
&\|\p_{y}^j \int_0^t\bar{\mathcal {U}}_{(0)}(t-\tau) \bar f_2(\tau,y)d\tau\|_{L^2}
\le C\|\p_{y}^{j} \int_0^t\bar{\mathcal {U}}_{(0)}(t-\tau) \p_y^2 (\f{\phi}{\gamma+\phi}m)d\tau\|_{L^2} \\
&\le C\int_\f{t}{2}^t(1+t-\tau)^{-\f{3}{2}-\f{j}{2}} \|\phi m\|_{L^1}d\tau
+ C\int_0^\f{t}{2}(1+t-\tau)^{-1-\f{j}{2}} \|\p_y (\phi\cdot m)\|_{L^1}d\tau \\
& \le C\int_\f{t}{2}^t(1+t-\tau)^{-\f{3}{2}-\f{j}{2}}\ve^2 (1+\tau)^{-1}d\tau
+ C\int_0^\f{t}{2}(1+t-\tau)^{-1-\f{j}{2}}(1+t)^{-\f{3}{2}}(\tau)d\tau \\
& \le C\ve^2 (1+t)^{-1}.\tag{5.4}
\end{align*}
In addition, one has
\begin{align*}
\bar f_1 -h_1&=\left(
\begin{array}{cccc} 0 \\ -\Div' (\Pi[\f{\gamma}{\gamma+\phi}m_1m'-w_1 w'])\\ -\Div' (\Pi[\f{\gamma}{\gamma+\phi}m'\otimes m' -w'\otimes w'])-\nabla ' \Pi[F(\phi)-\alpha \sigma^2]\end{array}
\right)\\
&
=\left(
\begin{array}{cccc}  0 \\ -\Div' (\Pi[\f{\gamma}{\gamma+\phi}m_1m'-\bar m _1 \bar m'])\\ -\Div' (\Pi[\f{\gamma}{\gamma+\phi}m'\otimes m' -\bar m'\otimes \bar m'])-\nabla ' \Pi[F(\phi)-\alpha {\bar \sigma}^2]\end{array}
\right)
\\
&\qquad +
\left(
\begin{array}{cccc} 0 \\ -\Div' ( \bar m _1 \bar m'-w_1 w')\\ -\Div' (\bar m'\otimes \bar m' -w'\otimes w')-\nabla ' (\alpha {\bar \sigma}^2-\alpha \sigma^2)\end{array}
\right)\\
&= R_1(t,y)+R_2(t,y).\tag{5.5}
\end{align*}
We now analyze the first term $R_1(t,y)$.

Since
$$
\f{\gamma}{\gamma+\phi} m'\otimes m' -\bar m'\otimes \bar m'=\f{1}{\gamma+\phi} (\gamma \tilde m'\otimes m'
 +\gamma \bar m'\otimes\tilde m'-\phi\bar m'\otimes\bar m')\eqno{(5.6)}
$$
and
\begin{align*}
&\|\gamma \tilde m'\otimes m'
+\gamma \bar m'\otimes \tilde m'\|_{L^1}\le C\|\tilde m' \|_2\|m'\|_2\le C(1+t)^{-\f{3}{2}}\ve^2,\tag{5.7}\\
&\|\f{1}{\gamma+\phi}\phi\bar m'\otimes\bar m'\|_{L^1}\le C\|\phi\|_{L^2}\|\bar m'\|_{L^\i}\|\p_y\bar m'\|_{L^2}
\le C(1+t)^{-\f{31}{12}}\ve^{\f{67}{24}},\tag{5.8}
\end{align*}
substituting (5.7)-(5.8) into (5.6) yields
$$
\|\f{\gamma}{\gamma+\phi} m'\otimes m' -\bar m'\otimes \bar m'\|_{L^1} \le C(1+t)^{-\f{3}{2}}\ve^2.\eqno{(5.9)}
$$

On the other hand, if we denote $\alpha_1(\phi)=\ds\f{1}{\g^2}\int_0^1(1-\th)^2p^{''}(1+\f{\th\phi}{\g})d\th $, then
$$
F(\phi)-\alpha\bar \phi^2=\alpha_1(\phi)(\tilde \phi \phi+\bar \phi \tilde \phi)+(\alpha_1(\phi)-\alpha)\bar \phi ^2.\eqno{(5.10)}
$$

Noting
$$
\alpha_1(\phi)-\alpha=\f{1}{\g^2}\int_0^1(1-\th)^2p^{''}(1+\f{\th\phi}{\g})d\th-\alpha=O(\phi),
$$
then we have
$$
\|F(\phi)-\alpha\bar \phi^2\|_{L^1}\le C(1+t)^{-\f{3}{2}}\ve^2.\eqno{(5.11)}
$$
Combining (5.9) with (5.11) yields
$$
\|\f{\gamma}{\gamma+\phi} m'\otimes m' -\bar m'\otimes \bar m'\|_{L^1}+\|\f{\gamma}{\gamma+\phi} m_1 m' -\bar m_1 \bar m'\|_{L^1}+
\|F(\phi)-\alpha\bar \phi^2\|_{L^1}\le C(1+t)^{-\f{3}{2}}\ve^2.
$$
From this and the expression of $R_1(t,y)$, we obtain for $j=0,1$,
\begin{align*}
&\|\p_y^j \int_{0}^t \bar {\mathcal {U}}_{(0)} (t-\tau )R_1(\tau,y)d\tau  \|_{L^2}\\
&\le C\int_{0}^t (1+t-\tau)^{-\f{j+2}{2}}(1+\tau)^{-\f{3}{2}}\ve^2d\tau \le C(1+t)^{-1}\ve^2.\tag{5.12}
\end{align*}
Finally, we estimate $R_2(t,y)$.
We see that for any fixed $\dl>0$
\begin{align*}
&\|w_1w' -\bar m_1 \bar m' \|_{L^1} \le \|w_1-\bar m_1\|_{L^2}\|\bar w'\|_{L^2} +\|w'-\bar m'\|_{L^2}\|\bar m_1\|_{L^2}
\le C(1+t)^{-\f{3}{2}}\ve N_1(t),\\
&\|\sigma^2 -\bar \phi^2 \|_{L^1} \le \|\sigma-\bar \phi \|_{L^2}\|\sigma+\bar \phi\|_{L^2}\le C(1+t)^{-\f{3}{2}}\ve N_1(t),\\
&\|w'\otimes w' -\bar m'\otimes  \bar m' \|_{L^1} \le \|w'-\bar m'\|_{L^2}\|\bar w'\|_{L^2} +\|w'-\bar m'\|_{L^2}\|\bar m'\|_{L^2} \le C (1+t)^{-\f{3}{2}+\delta}\ve N_1(t).
\end{align*}
This, together with the expression of $R_2(t,y)$, yields
\begin{align*}
&\|\p_y^j \int_{0}^t \bar {\mathcal {U}}_{(0)} (t-\tau ) R_2(\tau, y) d\tau  \|_{L^2}\\
&\le C N_1(t)\int_{\f{t}{2}}^t (1+t-\tau)^{-1-\f{j}{2}}(1+\tau)^{-\f{3}{2}+\dl}\ve d\tau +C N_1(t)\int_0^{\f{t}{2}}  (1+t-\tau)^{-1-\f{j}{2}}(1+\tau)^{-\f{3}{2}+\dl}\ve d\tau\\
& \le C(1+t)^{-1}\ve N_1(t).\tag{5.13}
\end{align*}
From (5.12)-(5.13), (5.5) and (5.4), we arrive at
$$
\|\bar u-\eta\|_{H^1} \le C(1+t)^{-1}\ve+ C(1+t)^{-1}\ve N_1(t),
$$
which derives for small $\ve >0$
$$
N_1(t) \le C\ve.
$$
Then we complete the proof of Lemma 5.1.\qquad $\square$

\vskip 0.2 true cm

{\bf Proof of Theorem 1.1.} From Remark 1.1, we know that  (1.3) has a global solution
$u(t,z)=(\phi,m)\in C([0, +\infty), H^4(\Bbb T \times \Bbb R^2))\cap
C^1([0, +\infty), H^2(\Bbb T \times \Bbb R^2))$. In addition, (1.4) comes from (4.13) directly,
and (1.5) comes from the interpolation between (1.4) and $\|u\|_{H^4}\le C\ve$, and (1.6)
is yielded by combining (4.14) and Lemma 5.1.\qquad $\square$

\vskip 0.4 true cm \centerline{\bf \S6. The proofs of Theorem 1.2 and Theorem 1.3}
\vskip 0.3 true cm

In this section, at first we focus on the proof of Theorem 1.2. To this end,
we will derive some decay properties of the solution to (1.3) for $u_0(z)\in H^4(\Bbb T^2 \times \Bbb R)$
as in $\S 4$.
Define
$$
\Pi [g](y)=\bar g(y) =\f{1}{(2\pi)^2}\int_{\Bbb T^2}g(x_1,x_2,y)dx_1dx_2.
$$
Then as in (4.3), the mean value ${\bar u}(t,y)$ of the solution $u(t,z)$ to (1.3) over $\Bbb T^2$ satisfies
a similar nonlinear partial differential system.
Let $\tilde u(t,z)=u(t,z)-\bar u(t,y)$ and
$$
M_1(t)=M_1[u](t)=\sup_{0\leq\tau\leq t}\{(1+\tau)^{\f{1}{4}}\|u(\tau,y)\|_{L^2(\Bbb T^2\times\Bbb R)}+(1+\tau)^\f{3}{4}\|\na u(\tau,y)\|_{L^2(\Bbb T^2\times\Bbb R)}\}.
$$

By the analogous proof in Lemma 4.2, we can obtain

{\bf Lemma 6.1. }{\it Assuming $u\in C([0, +\infty), H^4(\Bbb T^2\times\Bbb R)\cap L^1(\Bbb T^2\times\Bbb R))$,
$\|u\|_{H^4\cap L^1}\leq\ve$ and $M_1(t) \le 1$, then we have for small $\ve>0$

(i) $\|u\|_{L^\i}\leq C(1+t)^{-\f{1}{2}}M_1(t)+C\ve^\f{1}{2}(1+t)^{-\f{3}{8}}M_1^\f{1}{2}(t)$.

(ii) $\| \mathcal {G}(\phi,m)\|_{L^1}\leq C(1+t)^{-\f{1}{2}}M_1^2(t)$.

(iii) $\| g(\phi,m)\|_{L^1}\leq C(1+t)^{-\f{1}{2}}M_1^2(t)$.

(iv) $\| G\|_{H^1}\leq C(1+t)^{-\f{5}{4}}M_1^2(t)+C\ve^\f{1}{2}(1+t)^{-\f{3}{4}}M_1(t)$.

Where the definitions of $\mathcal {G}(\phi,m), g(\phi,m)$ and $ G$ have been given at the beginning
of $\S 4$}

\vskip 0.2 true cm

Based on Lemma 6.1, by applying Remark 3.1 and Lemma 4.1 (ii), as in the proofs of Lemma 4.4 and Lemma 4.5, we have

{\bf Lemma 6.2. }{\it Assuming $u\in C([0, +\infty), H^4(\Bbb T^2\times\Bbb R)\cap L^1(\Bbb T^2\times\Bbb R))$,
$\|u\|_{H^4\cap L^1}\leq\ve$ and $M_1(t) \le 1$, then
\begin{align*}
&\text{(i)}\quad M_1[\bar u](t)\leq C(\ve+\ve M_1(t)+M_1^2(t)).\\
&\text{(ii)}\quad  M_1[\tilde u](t)\leq C(\ve+\ve^\f{1}{2}M_1(t)+ M_1^{2}(t)).
\end{align*}
}

From Lemma 6.2, we can easily derive the following results similar to Proposition 4.6:

{\bf Proposition 6.3. }{\it For small $\ve>0$, if $u_0\in H^4(\Bbb T^2 \times \Bbb R)$
and $\|u_0\|_{H^4\cap L^1}\leq\ve$, then we have
$$ \|\p_z^k u(t,z)\|_{L^2} = O(t^{-\f{1}{4}-\f{k}{2}}),\ \ k=0,1,$$
$$ \|\p_z^k( u(t,z)-\bar u(t,y))\|_{L^2} = O(t^{-\f34}),\ \ k=0,1,$$
where $u$ is the solution of (1.3) and $\bar u$ is the mean value of $u$ over $\Bbb T^2$.}

In addition, as in Lemma 5.1, one has

{\bf Lemma 6.4. }{\it   If $u_0\in H^4(\Bbb T^2 \times \Bbb R)$ and $\|u_0\|_{H^4\cap L^1}\le \ve$, then for small $\ve>0$
and any fixed constant $\dl>0$
$$
\|\bar u -\eta\|_{H^1} \le C (1+t)^{-\f{3}{4}+\delta} \ve.
$$
}

{\bf Proof of Theorem 1.2.} Since we have established the crucial Proposition 6.3
and Lemma 6.4, the proof of Theorem 1.2 is completely analogous to that for
Theorem 1.2, we omit the details here.\qquad \qquad \quad $\square$

\vskip 0.2 true cm

Finally, we prove Theorem 1.3.

\vskip 0.2 true cm
{\bf Proof of Theorem 1.3.}  By (2.7), it is easy to see that the meanvalues $\bar u$ of $u(t,z)$
over $\Bbb T^3$ are constants, namely, $(\bar\phi,\bar m)\equiv(\bar\phi_0,\bar m_0)$. Denote $\tilde u=u-\bar u$,
then it follows from (1.3) that
\begin{equation}
\left\{
\begin{aligned}
&\partial_{t}\tilde\phi+\g\Div\tilde m=0,\\
&\p_t\tilde m+\biggl(-\ds\f{\g}{(\bar\phi+\g)^2}\bar m(\bar m\cdot\na)+\f{c}{\g}\na+\f{\nu_1\g}{(\bar\phi+\g)^2}\bar m\Dl+\f{\nu_2\g}{(\bar\phi+\g)^2}(\bar m\cdot\na)\na\biggr)\tilde\phi\\
&\qquad +\biggl(\ds\f{\g}{\bar\phi+\g}\bar m\Div+\f{\g}{\bar\phi+\g}\bar m\cdot\na-\f{\g}{\bar\phi+\g}(\nu_1\Dl+\nu_2\na\Div)\biggr)\tilde m\\
&\qquad \quad =\na[\tilde u^2g_1(\bar u,\tilde u)]+\na^2[\tilde u^2g_2(\bar u,\tilde u)],
\end{aligned}
\right.\tag{6.1}
\end{equation}
where $c=p'(\f{\bar\phi+\g}{\g})$, and $g_1$, $g_2$ are smooth functions on their argumants.

If we denote
\begin{align*}
&A=\\
&\left(
\begin{array}{cccc} 0&\g\Div\\-\f{\g}{(\bar\phi+\g)^2}\bar m(\bar m\cdot\na)+\f{c}{\g}\na+\f{\nu_1\g}{(\bar\phi+\g)^2}\bar m\Dl+\f{\nu_2\g}{(\bar\phi+\g)^2}(\bar m\cdot\na)\na&\f{\g}{\bar\phi+\g}\bigl(\bar m\Div+\bar m\na-(\nu_1\Dl+\nu_2\na\Div)\bigr)\end{array}
\right),
\end{align*}
then the linearized system of (6.1) can be written as
\begin{equation}
\left\{
\begin{aligned}
&\p_t v+Av=0,\\
&v(0,z)=\t u_0(z)=u_0(z)-\bar u_0.
\end{aligned}
\right.\tag{6.2}
\end{equation}
Taking the Fourier series expansion for the matrix $A$, we get
\begin{align*}
&\hat A(k)=\\
&\left(
\begin{array}{cccc} 0&i\g k^T\\-\f{i\g}{(\bar\phi+\g)^2}(\bar m\cdot k)\bar m+\f{ic}{\g}k-\f{\nu_1\g|k|^2}{(\bar\phi+\g)^2}\bar m-\f{\nu_2\g}{(\bar\phi+\g)^2}(\bar m\cdot k)k&\f{i\g}{\bar\phi+\g}\bar m(k^T+k)+\f{\g}{\bar\phi+\g}(\nu_1|k|^2+\nu_2kk^T)\end{array}
\right),
\end{align*}
where $k\in\Bbb Z^3$.

It is easy to know that the eigenvalues of $\hat A(k)$ are
$$\la_1=\la_2=\f{\g}{\bar\phi+\g}(\nu_1|k|^2+i\bar m\cdot k),\ \la_\pm=\f{\f{\g}{\bar\phi+\g}((\nu_1+\nu_2)|k|^2+2i\bar m\cdot k)\pm\sqrt{\f{\g^2(\nu_1+\nu_2)^2|k|^4}{(\bar\phi+\g)^2}-4c|k|^2}}{2}.$$
Consequently, for small $(\bar\phi_0,\bar u_0)$ and $k\neq0$, there exists a positive constant $a_0$ such that $\Re e\la_\kappa\geq a_0$
for $\kappa=1,2,\pm$, and we can obtain
$$
\|v(t,z)\|_{H^1}\leq Ce^{-a_0t}\|\tilde u_0\|_{H^1}.\eqno{(6.3)}
$$

Denote $M_2(t)= \sup \limits_{0\le \tau \le t } e^{a_0t}\|\tilde u\|_{H^1}$. Noting that $\|u\|_{H^4}\le \ve$
and further $\|\t u\|_{H^4}\le \ve$ hold, then we can arrive at
\begin{align*}
\|\na[\tilde u^2g_1(\bar u,\tilde u)]\|_{H^1}\le &C\|\tilde uD\tilde u\|_{L^2}+C\|(D\tilde u)^2\|_{L^2}+C\|\tilde uD^2\tilde u\|_{L^2}\\
\le &C\|\tilde u\|_{L^6}\|D\tilde u \|_{L^3}+C\|D\tilde u\|_{L^4}^2+C\|\tilde u\|_{L^6}\|D^2\tilde u \|_{L^3}\\
\le &C\ve^\f{1}{6}\|\tilde u\|_{H^1}^{\f{3}{2}}
\end{align*}
and
\begin{align*}
\|\na^2[\tilde u^2g_2(\bar u,\tilde u)]\|_{H^1}\le &C\|(D\tilde u)^2\|_{L^2}+C\|(D\tilde u)^3\|_{L^2}+C\|D\tilde uD^2\tilde u\|_{L^2}
+C\|\tilde uD^3\tilde u\|_{L^2}\\
\le &C\|D\tilde u\|_{L^4}^2+C\|D\tilde u\|_{L^6}^3+C\|D\tilde u\|_{L^3}\|D^2\tilde u \|_{L^6}+C\|\tilde u\|_{L^6}\|D^3\tilde u \|_{L^3}\\
\le & C\ve^\f{1}{2} \|\tilde u\|_{H^1}^\f{7}{6}.
\end{align*}
Hence we have from (6.1) and (6.3)
\begin{align*}
\|\t u\|_{H^1}\le &\|v\|_{H^1}+\int_0^t \|e^{-(t-\tau)A}\bigl(0,\na[\tilde u^2g_1(\bar u,\tilde u)]+\na^2[\tilde u^2g_2(\bar u,\tilde u)]
\bigr)(\tau, y)\|_{H^1}d\tau \\
&\le Ce^{-a_0t}\|u_0\|_{H^1}+\ve^\f{1}{6} \int_0^t e^{-a_0(t-\tau)}  \|\tilde u(\tau, y)\|_{H^1}^\f{7}{6}d\tau\\
&\le Ce^{-a_0t}\|u_0\|_{H^1}+\ve^\f{1}{6}e^{-a_0 t}M_2^\f{7}{6}(t)\int_0^t e^{-\f{1}{6}a_0\tau}d\tau \\
&\le Ce^{-a_0 t}\ve +\ve^\f{1}{6}e^{-a_0 t} M_2^\f{7}{6}(t)
\end{align*}
and further obtain
$$
M_2(t) \le C \epsilon+C\epsilon^\f{1}{6} M_2^\f{7}{6}(t).
$$
This derives that $M_2(t) \le C \ve$ and $\|\t u(t, y)\|_{H^1}$ decays exponentially with respect to the time $t$. Therefore,
we complete the proof of Theorem 1.3.\qquad\qquad \qquad \qquad \qquad \qquad \qquad  $\square$

\vskip 0.4 true cm

\end{document}